\documentclass[12pt]{article}
\usepackage{a4}
\usepackage{amssymb}

\begin{document}
\newtheorem{theorem}{Theorem}[section]
\newtheorem{lemma}[theorem]{Lemma}
\newtheorem{corollary}[theorem]{Corollary}
\newtheorem{definition}[theorem]{Definition}
\newtheorem{proposition}[theorem]{Proposition}
\newtheorem{defprop}[theorem]{Definition-Proposition}
\newtheorem{example}[theorem]{Example}
\newtheorem{remark}[theorem]{Remark}
\newcommand{\Proof}{\noindent{\bf Proof:} }
\catcode`\@=11
\@addtoreset{equation}{section}
\catcode`\@=12
\renewcommand{\theequation}{\arabic{section}.\arabic{equation}}

\def\sqr#1#2{{\vcenter{\vbox{\hrule height.#2pt\hbox{\vrule width.#2pt
  height#1pt \kern#1pt \vrule width.#2pt}\hrule height.#2pt}}}}
\def\square{\mathchoice\sqr64\sqr64\sqr{2.1}3\sqr{1.5}3}
\def\cove{{\cal A}}
\def\id{{\rm id}}
\def\infigamma{{\cal I}_\gamma^{\infty}}
\def\II{{\cal I}}
\def\FF{{\cal F}}
\def\EE{{\cal E}}
\def\GG{{\cal G}}
\def\al{\alpha}
\def\be{\beta}
\def\ga{\gamma}
\def\de{\delta}
\def\om{\omega}
\def\Om{\Omega}
\def\iy{\infty}
\def\intiy{\int_{-\ga\cdot\iy}^{\ga\cdot\iy}}
\def\pa{\partial}
\def\som{{\,{\rm s}\om}}
\def\RR{{\bf R}}
\def\CC{{\bf C}}
\def\ZZ{{\bf Z}}
\def\Zplus{\ZZ_{\ge0}}
\def\HH{{\cal H}}
\def\HD{\HH^{\rm D}}
\def\HS{\HH^{\rm S}}
\def\MM{{\cal M}}
\def\GG{{\cal G}}
\def\Re{{\rm Re}}
\def\Im{{\rm Im}}
\def\half{{1\over 2}}
\def\mPP{\medbreak}
\def\infigamma{\II_\ga^\iy}
\def\sLP{\smallbreak\noindent}
\def\mLP{\medbreak\noindent}
\def\qed{\hfill$\square$}
\def\m{{\rm m}}
\def\d{{\rm d}}
\def\ux{{\underline x}}
\def\uz{{\underline z}}
\def\upar{{\underline \partial}}
\def\xe{x_1^{e_1}\cdots\,x_n^{e_n}}
\def\xf{x_1^{f_1}\cdots\,x_n^{f_n}}
\def\de{\partial_n^{e_n}\cdots\,\partial_1^{e_1}}
\def\df{\partial_n^{f_n}\cdots\,\partial_1^{f_1}}

\def\uzed{{\underline z}}
\def\ut{{\underline t}}
\def\te{t_1^{e_1}\cdots\,t_n^{e_n}}
\def\com{{\cal E}}
\def\cove{{\cal A}}
\def\vece{{\cal B}}
\def\somE{{s{\cal E}}}

\title{On the $q$-convolution on the line}
\author{Giovanna Carnovale\thanks{Current email address: carnoval@mathy.jussieu.fr}\\
D\'epartement de math\'ematiques\\
Universit\'e Cergy-Pontoise\\
2, Avenue Adolphe Chauvin\\
Cergy-Pontoise Cedex, France\\ }

\date{December 14, 1999}
\maketitle

\begin{abstract}
I continue the investigation of a $q$-analogue of the convolution on
the line started in a joint work with Koornwinder and based on a
formal definition due to Kempf and Majid.
 Two different ways of  approximating functions by means of
the convolution and convolution of delta functions are introduced.  A new
family of functions  that forms an increasing chain of algebras depending on a
parameter $s>0$ is constructed. The value of the parameter for which
the mentioned algebras are  well behaved, commutative and unital is found.  In
particular a privileged algebra of functions belonging to the above
family  is shown to be  the quotient of an algebra studied in the
previous article modulo the kernel
of a $q$-analogue of the Fourier transform. This result has an analytic
interpretation in terms of analytic functions whose $q$-moments
have a particular behaviour. The same result  makes it possible  to extend
results on invertibility of the $q$-Fourier transform due to Koornwinder. A
few results on invertibility of functions with respect to the $q$-convolution
are also obtained and they are related  to 
solving certain simple linear  $q$-difference equations with  polynomial
coefficients. \end{abstract}

\section{Introduction}\label{intro}

In \cite{giotom} a $q$-analogue of the convolution on the line was defined
(inspired by results in \cite{KeMa}), a few algebras under this new convolution
were constructed and  commutativity of these algebras was investigated. \\
In particular, it was shown how commutativity of the $q$-convolution
strongly relies on whether a function can be uniquely determined  by its
$q$-moments.  The $q$-convolution defined therein had as a formal limit
the usual  convolution of $C$-valued functions on $\bf R$. The $q$-Fourier
transform studied in \cite{Koo},  involving $E_q$,  intertwines the
$q$-convolution product and ordinary product of functions.\\
The definition of the $q$-convolution was motivated by the results
in \cite{KeMa}, where Fourier transforms and convolution  product were defined
for braided covector algebras. However,  \cite{giotom} and the present paper
are  developed in a commutative setting and from an analytic point of view. 
On the other hand, since the braided line is  commutative as an algebra,
one could interpret many of our results as living in the braided
setting.  Remarks referring to the theory of braided groups and explaining
results from a braided theoretical  point of view appear in
this paper wherever the author thought they could help
to motivate definitions and results. They can easily be skipped by a
reader who is only interested in  classical $q$-analysis.\\
After recalling the basic definitions and the main results in \cite{giotom},
I prove properties of the $q$-convolution related to
$q$-integrability, approximation of functions by means of the $q$-convolution, 
and existence of zero divisors in the family of algebras studied in
\cite{giotom}.\\
I  construct a new algebra  of functions with respect to the $q$-convolution
that will behave more like an algebra of distributions than  like  an algebra
of functions. It  consists of an increasing chain of algebras depending on a
continuous parameter $s$.  The elements will all be given by the $q$-Gaussian
$e_{q^2}(-X^2)$ times an entire function $f$ whose coefficients have a
particular decreasing behaviour. In the terminology of \cite{ramis},  where
$q$ corresponds to our $q^{-1}$, this corresponds to the fact that  $f$ is a
function of $q$-exponential  growth of order $1$ and finite type (they
type is related to $s$). Showing that  on the algebras associated to a parameter  $s<1$ the
$q$-moment problem is determined, I  prove that those algebras are commutative
and algebraic domains (i.e. they have no zero-divisors). They are also unital
if $s\ge q^{\half}$. In particular, the existence of a unit is a peculiar
phenomenon of the $q$-case, showing that on this space the $q$-convolution may
be seen as a mid-way between convolution of functions and convolution of
distributions.  In particular, the unit can be expressed in terms of 
Jackson's $q$-Bessel function ${\cal J}^{(1)}_{\half}$ (see \cite{KoSw} and
references therein). Existence of a unit element and existence of zero
divisors were not treated extensively in \cite{giotom}, although they
implicitly appear in the examples  related to commutativity. It is easy to 
show that if an algebra under $q$-convolution has no zero divisors, then it
will be commutative, but the converse is not always true. \\
I also show that
the commutative, unital algebras corresponding to the value of the parameter
$s$ ranging in $[q^\half,\,1)$ must coincide. I provide a constructive
method that associates to power series with a good behaviour, a function in
this algebra having the given power series as generating series of its
$q$-moments. This will show that the algebra corresponding to
$s\in[q^{\half},1)$ is  isomorphic to the   the space of definition of a
formal version of  the $q$-Fourier transform quotiented by its kernel (as
homomorphism). The formal $q$-Fourier transform coincides with Koornwinder's
one on a particular subspace, namely, the ideal generated by the $q$-Gaussian
$e_{q^2}(-x^2)$ with respect to the $q$-convolution. I  use this results in
order to construct a new operator inverting the formal $q$-Fourier transform,
hence Koornwinder's $q$-Fourier transform. It coincides with the inversion
operator described by Koornwinder in \cite{Koo} on their common domain of
definition, but it is defined on a bigger space than the domain defined in
\cite{Koo}. One can use this new inversion formula  to obtain new relations
between bases of various spaces. Other inversion formulas of the $q$-Fourier
transform were obtained in the braided context and on different spaces in
\cite{KeMa} and in \cite{OR}.\\  Since the algebra I  constructed is unital, 
the question of invertibility with respect to the $q$-convolution araises.
Although the inverse does not necessarily belong to the algebra, I obtain a
few results  on invertibility in a somewhat extended algebra. I  show how
similarly to the classical case, inversion of  a function with respect to the 
$q$-convolution, can be interpreted in terms of  solving inhomogeneous
$q$-difference equations with constant coefficients (i.e. particular
$q$-difference equations with coefficients in ${\bf C}[x^{-1}]$). Those
equations can be transformed into $q$-difference equation with polynomial
coefficients, that are regular singular at $0$ but not at
infinity (see \cite{ramis} for the definition of regular singular), and  whose
characteristic equation has roots  $1,\, q,\,q^2,\,\ldots,\,q^{n-1}$  if the
equation has order $n$. Solutions of homogeneous $q$-difference equations
of a class including the above class  were already described by
Adams in \cite{adams2} and  methods for solving inhomogeneous
equations are given in  \cite{adams1}. See also \cite{adams3} for a survey on
what was known on $q$-difference equations in the thirties. A solution of  a
$q$-differential equation  with constant coefficients can therefore clearly be
found without using the $q$-convolution, but in particular cases
$q$-convolution can simplify the problem. Moreover, it can be used in order
to  determine the space of functions a solution can belong to and unicity of
the solution in a given space of functions. In particular, if $F$ is a
function of the form $e_{q^2}(-x^2)$ times a function of $q$-exponential
growth of order $1$ and  small enough finite type, I  can give conditions on
the $q$-differential operator  with constant coefficients $L$ under which the
solution $y$ of $Ly=F$  will be again of the same form as $F$.\\

\section{Definitions and
Notations}\label{defi} 

In this Section I recall the necessary background, the
notation  and the results of \cite{giotom}.\\
Let $q\in(0,1)$ be fixed. \\
Denote as usual
$(a;q)_k:=\prod_{j=0}^{k-1}(1-aq^{j})$,\quad
$(a;q)_{\iy}:=\lim_{k\to\iy}(a;q)_k\,$,\\
$[k]_{q}:={{1-q^{k}}\over{1-q}}\,$,\quad
$[k]_{q}!:={{(q;q)_k}\over{(1-q)^k}}\,$,\quad
$\bigl[{k\atop j}\bigr]_{q}:=
{{[k]_q!}\over{[j]_q![k-j]_q!}}={{(q;q)_k}\over{(q;q)_j(q;q)_{k-j}}}\,$.\\
For $q$-hypergeometric series the notation of
Gasper \& Rahman \cite{gasper} will be followed.\\
In particular, we  will need the functions
\begin{equation}e_q(x)={{1}\over{(x;q)_\iy}}=\sum_{k=0}^\infty{{x^k}\over{(q;q)_k}}\quad{\hbox {and}}\quad E_q(x)=
\sum_{k=0}^\infty{{q^{k\choose2}x^k}\over{(q;q)_k}}=(-x;q)_\iy\end{equation}
where the series expansion of $e_q(x)$ holds for $|x|<1$.\\ The $q$-derivative
of a function $f$  at $x\not=0$ is given by $(\pa
f)(x):={{f(x)-f(qx)}\over{(1-q)x}}$, and the  $q$-shift $Q$ of $f$ is given by
$(Qf)(x):=f(qx)$. For $\ga>0$, $L(\ga)$ denotes the  $q$-lattice
$\{\pm q^{k}\ga\,|\,k\in\ZZ\}$.\\
For a function $f$ on $L(\ga)$ the $q$-integral over $L(\ga)$
is denoted and defined by
\begin{equation}\label{q-int}
\int_\ga f=
\intiy f(t)\,d_qt:=
(1-q)\sum_{k=-\iy}^\iy\sum_{\epsilon=\pm1}q^{k}\ga\,f(\epsilon q^k\ga),
\end{equation}
provided the summation absolutely converges.\\
A function $f\colon x\mapsto f(x)$ may also be denoted as $f(X)$.
This will be useful for functions like
$fX^e\colon x\mapsto f(x)x^e$ and
$e_{q^2}(-X^2)\colon x\mapsto e_{q^2}(-x^2)$.\\
For $\ga>0$, $\II_\ga$ denotes the space of absolutely $q$-integrable
functions on $L(\ga)$, and 
$\II_\ga^\iy$ denotes the subspace of functions $f\in\II_\ga$
such that  $f X^e$ is in $\II_\ga$ for every $e\in\Zplus$.
For $f\in\infigamma$ define the {\em moments}, respectively
{\em strict moments} of $f$ by:
\begin{equation}\label{momenteq}
\mu_{e,\ga}(f):=q^{{e^2+e}\over2}\intiy f(x)x^e\,d_qx,
\quad
\nu_{e,\ga}(f):=q^{{e^2+e}\over2}\intiy |f(x)x^e|\,d_qx.
\end{equation}
\sLP
For a real number $\al>0$, $\II_{\ga,\al}^\om$ denotes the space of
functions {\em of left type} $\al$ on $L(\ga)$ consisting of
all $f\in \infigamma$
such that, for some $b>0$,
$|\mu_{e,\ga}(f)|=O(q^{{\al e^2}\over2}b^e)$
as $e\to\iy$ and  $\II_{\ga,\al}^\som$ denotes the space of
functions {\em of strict left type} $\al$ on $L(\ga)$ consisting of
all $f\in \infigamma$
such that, for some $b>0$,
$\nu_{e,\ga}(f)=O(q^{{\al e^2}\over2}b^e)$
as $e\to\iy$.\\ 
The space $\II_\ga^\om$ (resp. $\II_\ga^\som$) denotes the union of all
$\II_{\ga,\al}^\om$ (resp. $\II_{\ga,\al}^\som$).  The space
$\II^\om_{\ga,>\al}$ (resp. $\II^\som_{\ga,>\al}$) denotes the union of 
all $\II^\om_{\ga,\beta}$ (resp. $\II^\som_{\ga,\beta}$) for
$\beta>\alpha$.\\  By $\HD$  (resp $\HS$) I  denote  the space of all functions
which are holomorphic on some disk (resp. strip) centered in 0 (resp. 
$\RR$).  Observe that if $f\in\HD$ then $\pa f$ is defined also at $x=0$. $\cal
E$ will denote the space of entire analytic functions.\\

Every time I  write $\HS$ or $\HD$ followed by one of the spaces denoted
by $\II_\ga$ with some index  (i.e. $\II_\ga$,
$\infigamma$  or the spaces of strict left functions or left functions on
$L(\ga)$) I  mean the space of functions contained in some sort of
intersection.  For instance, $\HD\infigamma$ will
denote the space of the functions on $L(\ga)$ belonging to $\infigamma$ which
coincide within some disk centered in 0 with the restriction of a (necessarily
unique) holomorphic function on that disk, with the assumption that the $\{\pm
q^k\gamma\,|\,k\in\Zplus\}$ is contained in the disk.  This assumption is not
restrictive since $\int_\ga f= \int_{q^k\ga}f$.

Note that  a function in $\HD\infigamma$ cannot be uniquely determined by its
restriction on a disk, so that the data of the values of the function on
$L(\ga)$  outside the disk should always be added. 
We recall that if $f$ is of left type $\alpha$ or of strict left type $\al$,
then every polynomial times $f$ is again so.  

\begin{definition}\label{convodefi}
Let $f\in\infigamma$ and let $g$ be
a function defined on some subset of $\CC$.
Then the {\em $q$-convolution product} $f*_\ga g$ is the function
given by
\begin{equation}\label{convoeq}
(f*_{\ga}g)(x):=\sum_{e=0}^\iy{{(-1)^e
\mu_{e,\ga}(f)}\over{[e]_{q}!}}(\pa^eg)(x)
\end{equation}
for $x\in\CC$ such that the $q$-derivatives $(\pa^eg)(x)$ are well-defined
for all $e\in\Zplus$ and the
sum on the right converges absolutely.
\end{definition}
\medbreak
\noindent For an $f\in\II^\iy_\ga$ we will denote by $\mu_\ga(f)$ and
$\nu_\ga(f)$ the series $\sum_{k=0}^\iy{{\mu_{k,\ga}(f)t^k}\over{[k]_q!}}$ and
$\sum_{k=0}^\iy{{\nu_{k,\ga}(f)t^k}\over{[k]_q!}}$  respectively. In
particualr, in \cite{giotom} it was proved that for $f\in\HD\II^\om_{\ga'}$ and
$g\in\HD\II_{\ga}^\om$
\begin{equation}\label{generating}\mu_\ga(f*_{\ga'}g)=\mu_{\ga'}(f)\mu_{\ga}(g)=\mu_{\ga'}(g*_\ga f)\end{equation} 
Let $\II_\ga^{\cal E}$  and $\II_\ga^{s{\cal E}}$denote the space of functions
in $\II^\iy_\ga$ for which $\mu_{\ga}(f)\in{\cal E}$ and 
$\nu_{\ga}(f)\in{\cal E}$  respectively.  $\HD\II^{\cal
E}_\ga:=\HD\II^{\iy}_\ga\cap\II^{\cal E}_\ga$ and $\HS\II^{\cal
E}_\ga:=\HD\II^{\iy}_\ga\cap\II^{\cal E}_\ga$ and similarly, $\HD\II^{s{\cal
E}}_\ga:=\HD\II^{\iy}_\ga\cap\II^{s{\cal E}}_\ga$ and $\HS\II^{s{\cal
E}}_\ga:=\HD\II^{\iy}_\ga\cap\II^{s{\cal E}}_\ga$. One can check that formula
(\ref{generating}) still holds for  $f\in\HD\II^{\cal E}_{\ga'}$ and
$g\in\HD\II_{\ga}^{\cal E}$.\\ 
In \cite{giotom} the following results are
achieved:    
 \begin{proposition}\label{results} With the notation
just introduced: \begin{enumerate}
\item The class $\HD\II_\ga^\om$
is an algebra (not necessary unital) with respect to $*_{\ga}$.
Its subclass $\HS\II_\ga^\om$ is also an algebra (not necessary unital)
and it is a left ideal of $\HD\II_\ga^\om$.
\item The class $\HD\II_\ga^\som$  is a subalgebra 
of $\HD\II_\ga^\om$. Its subclass $\HS\II_\ga^\som$ is a
left ideal of $\HS\II_\ga^\om$, $\HD\II_\ga^\om$ and $\HD\II_\ga^\som$.
\item The classes $\HS\II_{\ga,>c}^\som$ (for $c\in[0,1)$)
and $\HS\II_{\ga,c}^\som$ (for $c\in(0,1]$) are 
left ideals of $\HS\II_{\ga}^\som$ and of $\HS\II_{\ga}^\om$.
Similar properties hold for $\HD\II_{\ga}^\som$. 
\item Let $f,g\in\HD\II_\ga^\om$, $h\in\HD$. Then
$(f*_{\ga}g*_{\ga}h)(x)=
(g*_{\ga}f*_{\ga}h)(x)$ for every $x$ where
the product is defined. In particular, 
for every pair of ideals $I\subset J$ with
$I,\,J\in\{\HS\II_\ga^\som,\,\HS\II_\ga^\om,\,
\HD\II_\ga^\om,\,\HD\II_\ga^\som\,\}$,
$I$ is a left module over $J/[J,\,J]_*$, where $[J,\,J]_*$
denotes the commutator ideal.
\item A subalgebra $A$ of $\HD\II^\om_\ga$ under $q$-convolution is commutative
if  every function in the commutator $[A,\,A]_*$ is determined by its
$q$-moments. In particular, the $q$-moment problem is determined 
on $\HS\II_{\ga,>\half}^\som$ so that  $\HS\II_{\ga,>c}^\som$ is a commutative
algebra for every $c\in[1/2,1)$. \hfill$\square$ \end{enumerate}
\end{proposition} 
Statements {\sl 1}, {\sl 2} and  {\sl 4} and the first part of statement {\sl
5} still hold if we replace $\HD\II_\ga^\om$, $\HD\II_\ga^\som$,
$\HS\II_\ga^\om$ and  $\HS\II_\ga^\som$ by $\HD\II_\ga^{\cal E}$,
$\HD\II_\ga^{s{\cal E}}$  $\HS\II_\ga^{\cal E}$  and $\HS\II_\ga^{s{\cal E}}$
respectively.  However, there is  clearly no $\cal E$-counterpart  for 
$\HS\II_{\ga,>\half}^\som$.  
\begin{proposition}\label{results2}
Let $(\FF_\ga f)(y):=\int_{-\ga.\iy}^{\ga.\iy} E_q(-iqxy)\,f(x)\,d_qx$ and
let \\
$(\tilde\FF_\ga f)(y):=
\sum_{k=0}^\iy\mu_{k,\ga}(f)\,{(-iy)^k\over(q;q)_k}=\mu_\ga(f)(-iy(1-q)^{-1})$.  
\begin{enumerate}
\item If
$f\in\II_\ga^\om$ then $\tilde\FF_\ga f$ is well-defined and it is an entire
analytic function. If moreover $f\in\II_\ga^\som$ then $\FF_\ga f$ is also
well-defined and $\FF_\ga f=\tilde\FF_\ga f$.
\item Let $f\in\HD\II_\ga^\om$, $g\in\HD\II_{\ga'}^\om$.
Then $f*_\ga g\in\HD\II_{\ga'}^\om$ and
$\tilde \FF_\ga(f*_\ga g)=(\tilde \FF_\ga f)(\tilde\FF_{\ga'} g)$, so
$\tilde\FF_\ga$ is an algebra homomorphism from  the algebra $\HD\II^\om_\ga$
with convolution product to the algebra $\cal E$ with ordinary product. Its
kernel is given by functions for which all $q$-moments are zero.
\item  $\tilde\FF_\ga$ is an injective algebra homomorphism on  all subspaces
of $\HD\II^\om_\ga$ on which the $q$-moment problem is determined. In
particular, it is injective on  $\HS\II^\som_{\ga,>\half}$, where it coincides
with $\FF_\ga$.
\qed \end{enumerate} \end{proposition}
Again, statements involving $\II^\om_\ga$ and $\II^{\som}_\ga$ still hold if we
replace the upper index  $\om$ by ${\cal E}$ and $s{\cal E}$. This will allow
us to extend the domain of the  formal $q$-Fourier transform $\tilde\FF_\ga$ 
and its inverse. 
I
will frequently make use of the following formulas, deduced from formulas
(9.8) and  (9.14) in \cite{Koo} (for $k\in\Zplus$): 
\begin{eqnarray}
\label{eq-moment}
\mu_{2k,\ga}(e_{q^2}(-X^2))=c_q(\ga)\,q^{k^2-k}\,(q;q^2)_k
\quad\quad\mu_{2k+1,\ga}(e_{q^2}(-X^2))=0\\
\label{Eq-moment}
\mu_{2k,1}(E_{q^2}(-q^2X^2))=b_q\,q^{2k^2+k}\,(q;q^2)_k\quad\quad\mu_{2k+1,1}(E_{q^2}(-q^2X^2))=0
\end{eqnarray}

where $b_q=\int_{-1\cdot\iy}^{1\cdot\iy}
E_{q^2}(-q^2x^2)\,d_qx$ and $c_q(\ga)=\intiy e_{q^2}(-x^2)\,d_qx$ 
were given \cite{Koo}.

\section{Properties of
the convolution}\label{properties}

In this Section I  describe useful properties of the convolution that 
were not discussed  in \cite{giotom}. In particular I  investigate  the 
behaviour of the convolution product with respect to integration and
approximation of  functions. 

\begin{definition} \label{antipode} Let $S$ be the linear map from $\HD$ to
$\cal E$ mapping $\sum_r a_rx^r$ to $\sum_r(-1)^r q^{r\choose 2} a_r x^r$.
Let  $\varepsilon$ be the map from $\HD$ to $\CC$  evaluating
a function at $0$.\end{definition} 

\begin{remark}\label{braidedline} \rm A motivation for the definition of
$S$ and $\varepsilon$ comes from the braided Hopf algebra structure that the
algebra of power series $\CC[[x]]$ has (see \cite{Maj1}, \cite{Majiprimo} and
\cite{Koor}).  The braiding isomorphism for this braided Hopf algebra is 
$\Psi(x^k\otimes x^l):=q^{kl} x^l\otimes x^k$ and the  comultiplication is 
$\Delta(x^k):=\sum_{j=0}^k\;\bigg[{k\atop j}\biggr]_q x^{k-j}\otimes x^j$.
In this setting, the operators just defined
are  the counit $\varepsilon$  given by
$\varepsilon(x^k):=\delta_{k,0}$ and the braided antipode $S$ mapping $x$ to
$-x$ and extended as a braided-antimultiplicative map.  This braided Hopf
algebra $\cal A$ is called  the {\em braided line}  and it is the simplest
example of  a braided covector algebra.\hfill$\spadesuit$
\end{remark}
\begin{remark}\rm 
The operator $S$ followed by the map $f(x)\to f(-x)$ is
Changgui Zhang's formal $q$-Borel transform (see \cite{Zhang}, where
$q>1$). In his article, Zhang used the $q$-Borel transform in order to
define summability of certain divergent power series (the so-called $q$-Gevrey
series originally introduced in \cite{bezi}) with a particular growth.
\hfill$\spadesuit$ \end{remark}
Recall that $Q$ is the $q$-shift operator such that $(Qf)(x)=f(qx)$.\\
\begin{proposition}\label{gcap} Let $f\in\II^\somE_{\ga,\al}$ and let
$g\in\HD$. Then for  every $p\in\ZZ$, $f\,Q^pSg\in\II^\iy_\ga$.
If $f,\,g$ and $h\in\HS\II_{\ga,>\half}^\som$   there holds 
\begin{eqnarray}\int_{\gamma}f\,(Q(Sg))&=&\int_{\gamma}(Q
(Sf))\,g\\
\label{property}\int_{\gamma}(f*_{\gamma}g)\,Q(S(h))&=&\int_{\gamma}(g*_{\gamma}h)\,Q(S(f)).
\label{classical}\end{eqnarray} so in particular 
the  restriction of $SQ=QS$  to $\HS\II_{\ga,>\half}^\som$ is symmetric with respect to the
$q$-integral. \end{proposition} {\bf Proof:}  Let $g\in\HD$ be such that 
$g(x)=\sum_{l=0}^\infty a_lx^l$ for $|x|<\rho$.  Since $Sg\in{\cal E}$,  for
any fixed $k\in\ZZ_{\ge0}$ and $p\in\ZZ$:
\begin{eqnarray*}&&\int_\gamma|X^k(Q^pSg)f|\le \int_\gamma\sum_{l=0}^\infty
q^{{{l^2-l}\over2}+lp}|a_l||f X^{l+k}|\\ &&\le C
q^{{\half  k^2}-kp}(\rho')^{-k}\sum_{l=0}^\infty
\bigl((\rho')q^{-\half l}q^{-k+p}\bigr)^{l+k}\nu_{l+k,\ga}(f)<\infty
\end{eqnarray*}
for some constants $C$ and $\rho'>\rho$ and the first statement is proved.\\
Let now  $f$ and $g$  in $\HS\II_{\ga,>\half}^\som$. By statement {\sl
5} of Proposition \ref{results} one has 
$\varepsilon(f*_{\gamma}g)=\varepsilon(g*_{\gamma}f)$. By statement {\sl 2} in
Proposition \ref{results}, $f*_{\gamma}g$ and $g*_{\gamma}f$ can be written as
power series for $|x|$ sufficiently small. In particular 
$$f*_{\gamma}g(x)=\sum_{r=0}^\infty\Biggl[\sum_{k=0}^\infty(-1)^k
\mu_{k,\gamma}(f)a_{r+k}
\Bigl[ {{r+k}\atop{k}}\Bigr]_q\Biggr]x^r$$ so that
$\varepsilon(f*_{\gamma}g)=\sum_{k=0}^\infty(-1)^kq^{{k^2+k}\over2}a_{k}\int_\ga f X^k$. Since $f$ is 
of strict left type we may interchange integration
and summation by dominated convergence. Hence
\begin{equation}\int_{\gamma}f\,(Q(Sg))=\varepsilon(f*_{\gamma}g)=\varepsilon(g*_{\gamma}f)
=\int_{\gamma}(Q(Sf))\,g\label{scambio}\end{equation} 
Formula (\ref{property}) follows by (\ref{scambio}) together with
associativity of the $q$-convolution because 
$$\int_{\gamma}\left(f*_{\gamma}g\right)\,(Q(Sh))=\varepsilon((f*_{\gamma}g)*_{\gamma}h)= \varepsilon(f*_{\gamma}(g*_{\gamma}h))=\int_{\gamma}(g*_{\gamma}h)\,(Q(Sf))$$
\hfill$\square$\smallbreak
\noindent Observe that formula (\ref{property})  
resembles the classical property that
\begin{equation}\int_{-\infty}^{\infty}(f*g)(x)h(-x)dx=\int_{-\infty}^{\infty}(g*h)(x)f(-x)dx.\end{equation}
An important classical property of the convolution is
that it makes it possible to approximate functions by  sequences of functions
with a nice behaviour. I  introduce here a method to approximate functions by
means of the $q$-analogue of the convolution. 
\begin{proposition}\label{approx} Let $f\in \II_\ga^{\cal E}$ be such that 
$\int_\ga f=1$. Define the  sequence $f_k(x):=q^{-k}Q^{-k}f\in \II_\ga^{\cal
E}$ for $k\in\Zplus$.  Let $g$ be a function defined on a domain $\Omega$
together with all its $q$-derivatives. If for every $x\in\Omega$ there exists
an $R_x>0$ such that $\bigl|\partial^kg(x)\bigr|=O(R_x^k)$ as $k\to\iy$,
then  $\lim_{l\to\infty}f_l*_\ga g(x)=g(x)$ pointwise for every
$x\in\Omega$. If the majorization of $|\partial^kg(x)\bigr|$ is uniform, then
the limit is uniform. In particular, this holds on a disk centered at zero  if
$g\in\HD$. \end{proposition}  {\bf Proof :}  Since
$\mu_{l,\ga}(f_k)=q^{kl}\mu_{l,\ga}(f)$, for every fixed $x\in\Omega$, the
product $f_l*_{\ga}g$  at $x$ is an absolutely convergent  sum, so that 
$$\biggl|f_l*_{\ga}g(x)-g(x)\biggr|=\Biggl|\sum_{r=1}^\infty
{{(-1)^{r}q^{lr}\mu_{r,\ga}(f)}\over{[r]_q!}}\pa^rg(x)\Biggr|\le
Cq^l\sum_{r=1}^\iy {{|\mu_{r,\ga}(f)|R_x^r}\over{[r]_q!}}$$
for some constant $C>0$.
Hence for every fixed $x\in\Omega$,  $\bigl|f_r*_\ga
g(x)-g(x)\bigr|\to 0$ as $l\to\iy$.\\
If for some $R>0$ 
one has $\bigl|\partial^kg(x)\bigr|=O(R^k)$ as $k\to\iy$ for every
$x\in\Omega$, then clearly $\bigl|\bigl|f_l*_\ga
g-g\bigr|\bigr|_\Omega\to 0$ as $l\to\iy$.
\hfill$\square$\smallbreak

One  would prefer   to approximate
functions by convolution from the right because $f*_\ga g$ inherits the
properties of $g$, not those of $f$. However one cannot have ``good''
convergence in this case in general, unless $f$ has good properties too.
Moreover,  it has been shown in Theorem 6.5 and Example 6.7 in \cite{giotom}
that there are nonzero functions in $\HD\II_\ga^\som$ and  in
$\HD\II_{\ga,\al}^\om$ for every $\al>0$ whose $q$-moments are all zero. Those
functions could never be approximated using $q$-convolution from the right. In
fact, the functions that could  nontrivially be approximated using
$q$-convolution from the right  by a  sequence of ``good'' functions are those
for which the $q$-moment problem is determined, i.e. the functions that commute
with  ``good'' ones. In that case, approximation from the left and from the
right coincide.

\begin{example}\rm Let $\ga\in(0,\,1)$  and
$f_\ga={{e_{q^2}(-X^2)}\over{c_q(\ga)}}$, with $c_q(\ga)$ is as in formula
(\ref{eq-moment}) and let $f_{k,\ga}=q^{-k}Q ^{-k}f_\ga$ defined on
$|\Im(x)|<q^k$. For any $g\in\HD$,  $g(x)=\sum_lc_lx^l$  on a
neighbourhood of $0$ and the sequence of products:
 \begin{eqnarray*}&&(f_{k,\gamma}*_{\gamma}g)(x)=
\sum_{l=0}^\infty(-i)^lc_l(q;q)_l
\sum_{e=0}^{[{l\over2}]}{{(-1)^eq^{e(e-1)}q^{2e(k+1)}(ix)^{l-2e}}\over{(q;q)_{l-2e}(q^2;q^2)_e}}=\\
&&\sum_{l=0}^\infty(-i)^lq^{l(k+1)}c_lh_{l}(iq^{-(k+1)}x;q)
\end{eqnarray*} (where the $h_l(x;q)$'s are the discrete $q$-Hermite $I$
polynomials, see \cite{KoSw}) converges to $g$ uniformly on a disk centered at
$0$ for $k\to\iy$. \\  
If $\gamma=1$ and 
$F={{E_{q^2}(-q^2X^2)}\over{b_q}}$ with $b_q$ as in formula
(\ref{Eq-moment}) and if  $F_{k}:=q^{-k}Q^{-k}F$ for $k\ge0$, the sequence of
products 
\begin{eqnarray*}&&(F_{k}*_1g)(x)=\sum_{l=0}^\infty(-i)^lq^{(k+l)l}c_l{\tilde
h}_l(iq^{-(k+l)}x;q)
\end{eqnarray*} (where the ${\tilde h}_l(x;q)$'s are the
discrete $q$-Hermite $II$  polynomials, see \cite{KoSw}) converges uniformly
to $g$ on a disk for $k\to\iy$.
\hfill{$\spadesuit$} \end{example}

We shall see now  how the $q$-convolution product is related
to the ordinary product of functions. One can prove that for  two
functions $f$ and $g$, $$\partial^n(fg)=\sum_{k=0}^n\biggl[{n\atop
k}\biggr]_q\bigl(Q^k\partial^{n-k}f\bigr) \,\partial^kg$$ (see  also T.
Koornwinder's informal note \cite{informal}). Then, for $f$ and $g$ in $\HD$, 
and $h\in\II_\ga^{\cal E}$, there holds \begin{eqnarray*}&& \bigl(h*_\gamma
(fg)\bigr)=\sum_{e=0}^\infty{{(-1)^e\mu_{e,\gamma}(h)}\over{[e]_q!}}
\sum_{k=0}^e\biggl[{e\atop k}\biggr]_q\bigl(Q^k\partial^{e-k}f\bigr)
\,\partial^kg\\
&&=\sum_{k=0}^\infty{{(-1)^kq^{{k^2+k}\over2}\partial^kg}\over{[k]_q!}}
\sum_{s=0}^\infty{{(-1)^{s}\mu_{s,\gamma}(X^k
h)}\over{[s]_q!}}\partial^s(Q^kf)\\
&&=\sum_{k=0}^\infty{{(-1)^kq^{{k^2+k}\over2}}\over{[k]_q!}}
\biggl(h\,X^k*_\ga Q^kf\biggr)\partial^k g.  \end{eqnarray*}
In particular, for $g=X$ this implies 
\begin{equation}X(h*_\ga f)=q(hX*_\ga
Qf)+h*_\ga(fX)\label{leibniz}\end{equation} i.e. multiplication by $X$ obeys a
sort of Leibniz rule. \\

The end of  this Section is devoted to a few remarks about functions whose product
is zero.
\begin{lemma} \label{domain}Let $f\in\HD\II_\ga^{\cal E}$ and 
$g\in\HD\II_{\ga'}^{\cal E}$. If $f*_\ga g=0$, then  $\mu_{\ga}(f)=0$ 
 and/or $\mu_{\ga'}(g)=0$.  Hence, a
convolution subalgebra  $A$ of $\HD\II_\ga^{\cal E}$ has no zero divisors if
and only if the $q$-moment problem is determined on $A$.  In particular,
$\HS\II^\som_{\ga,>\half}$ is an algebraic domain. \end{lemma} \Proof  This is
a trivial consequence of formula (\ref{generating}). The last statement
follows by Lemma 6.1 in \cite{giotom}.\hfill$\square$\smallbreak

\begin{corollary}\label{domecommu}Let $f\in\HD\II_\ga^{\cal E}$ and 
$g\in\HD\II_{\ga'}^{\cal E}$. If $f*_\ga g=0$, then either  $f*_\ga h=0$ for
every $h$ or $g*_{\ga'}h=0$ for every $h$ for which the product is defined. In
particular, either  $f$ or
$g$ is nilpotent.\end{corollary}
\Proof  If $\mu_{\ga}(f)\equiv0$  then $f*_\ga h=0$ for
every $h$.\hfill$\square$\smallbreak 

\begin{remark}\label{domicommu}\rm A consequence of the above results is
that  if a subalgebra $\cal N$ of  $\HD\II_\ga^\com$ has no zero divisors (i.e.
if the $q$-moment problem is determined on $\cal N$), then $\cal N$ is
commutative since by Corollary 5.11 in \cite{giotom} $(f*_\ga g-g*_\ga f)*_\ga
h=0$. The converse does not necessarily hold. Indeed, for $\ga<1$ consider the
functions $e_{q^2}(-q^2X^2)$ and $e_q(iX)$. By formula (8.21) in
\cite{Koo} for $t=q^{-1}$ one has $\mu_{\ga}(e_q(iX))\equiv0$. Hence,
 $e_q(iX)*_\ga f=0$ for every $f\in\HD$.
By 
formulas (\ref{eq-moment}) 
$\,\mu_{2r+1,\ga} (e_{q^2}(-q^2X^2))=0$ and 
$\mu_{2r,\ga}(e_{q^2}(-q^2X^2))=c_q(\ga) q^{-2}(q;q^2)_rq^{r^2-r}$ (see
also Section 9 in \cite{Koo}). Since  $\pa^k
e_q(ix)=i^k(1-q)^{-k}e_q(ix)$  one has
\begin{equation}e_{q^2}(-X^2)*_\ga e_q(iX)=q^{-2}
c_q(\ga)\sum_{r=0}^\infty{{(-1)^rq^{r^2-r}}\over{(q^2;q^2)_r}}e_q(iX)=
E_{q^2}(-1)e_q(iX)=0\end{equation}  Hence the subalgebra generated by
$e_{q^2}(-q^2X^2)$ and $e_q(iX)$ is commutative although it has zero
divisors.    \hfill$\spadesuit$\end{remark}

\section{Discrete delta functions}\label{distributions}

This Section is devoted to the study of the $q$-convolution for  discrete delta
functions, and is based on ideas of T. Koornwinder.\\
 On the $q$-lattice
$L(\gamma)$ define the discrete delta functions $$\delta_{\epsilon\gamma
q^p}(\eta \gamma q^l):=\delta_{\epsilon,\eta}\delta_{l,p}$$ for any
$\epsilon,\eta\in\{\pm 1\}$ and any $l$ and $p\in{\bf Z}$.  By Ryde's
formula (see \cite{ryde})  for iterated $q$-differentation
at $x\not=0$ 
$$(\partial^nf)(x)=(1-q)^{-n}x^{-n}\sum_{k=0}^n(-1)^k\biggl[{n\atop
k}\biggr]_{q} q^{-k(n-k)}q^{-{{k(k-1)}\over2}}Q^kf(x)$$ and since
$Q\delta_{\eta\gamma q^{t}}=\delta_{\eta\gamma q^{t-1}}$ one sees that 
$\partial^n\delta_{\eta\gamma q^{t}}$ is a linear combination of  discrete
$\delta$ functions  $\delta_{\eta\gamma q^{s}}$ for $t-n\le s\le t$. In
particular, for $t-n\le s\le t$ and $\epsilon\in\{\pm1\}$, 
$$(\partial^n\delta_{\eta\gamma q^{t}})(\epsilon \gamma
q^{s})=(-1)^{t-s}\delta_{\epsilon,\eta}{{\eta^n\gamma^{-n}q^{-sn}}\over{(1-q)^n}} \biggl[{n\atop{t-s}}\biggr]_{q}q^{-{{(t-s)(t-s-1)}\over2}
-(t-s)(n-t+s)}.$$
Moreover, $$\mu_{k,\ga}(\delta_{\eta\gamma
q^{n}})=(1-q)\eta^kq^{\half(k^2+k)}\gamma^{k+1}q^{n(k+1)}\quad{\hbox{and}}\quad\nu_{k,\ga}(\delta_{\eta\gamma
q^{n}})=|\mu_{k,\ga}(\delta_{\eta\gamma
q^{n}})|$$
hence the convolution
product of discrete delta functions along $L(\gamma)$ is well-defined by
Proposition \ref{results}. One computes, for $l\le s$ and $\theta\in\{\pm1\}$: 
\begin{eqnarray*}&&(\delta_{\epsilon\gamma q^{t}} *_{\gamma}
\delta_{\eta\gamma q^{s}})(\theta\gamma q^{l})=\sum_{k=0}^\infty {{(-1)^k
q^{{k^2+k}\over2}(1-q)\gamma^{k+1}\epsilon^k q^{t(k+1)}}
\over{[k]_{q}!}}\partial^k(\delta_{\eta\gamma q^{s}}(\theta\gamma q^{l}))\\
&&=\sum_{k=s-l}^\infty {{(-1)^k q^{{k^2+k}\over2}(1-q)\gamma^{k+1}\epsilon^k
q^{t(k+1)}}\over{[k]_{q}!}} \times\\
&&\times{{(-1)^{s-l}\delta_{\theta,\eta}
\eta^k\gamma^{-k}q^{-lk}}\over{(1-q)^k}}
\biggl[{k\atop{s-l}}\biggr]_{q} q^{-{{(s-l)(s-l-1)}\over2}-(s-l)(k-s+l)}\\
&&={{\gamma(1-q)\delta_{\theta,\eta}q^{(t+s-l)+(s-l)(t-l)}(\eta\epsilon)^{s-l}}\over{(q;q)_{s-l}}}
(\eta\epsilon q^{(t-l+1)};q)_{\infty}.\\
\end{eqnarray*}
Hence, if  $\eta=\epsilon$ one has
\begin{eqnarray*}(\delta_{\eta\gamma
q^{t}} *_{\gamma} \delta_{\eta\gamma q^{s}})(\theta\gamma
q^{l})={{\gamma(1-q)\delta_{\theta,\eta}q^{(t+s-l)+(s-l)(t-l)}}\over{(q;q)_{s-l}}} (q^{(t-l+1)};q)_{\infty}
\end{eqnarray*} which is zero for  $l>t$  so that the product is a linear
combination of  discrete delta functions with support $\eta\gamma q^{l}$ for 
$l\le \min(s,\,t)$. For $l\le t$  the product  evaluated at $\theta\gamma
q^{l}$ is equal to \begin{eqnarray*}
&&(\delta_{\eta\gamma q^{t}} *_{\gamma} \delta_{\eta\gamma
q^{s}})(\theta\gamma
q^{l})={{\gamma(1-q)\delta_{\theta,\eta}q^{(t+s-l)+(s-l)(t-l)}}\over{(q;q)_{s-l}(q;q)_{t-l}}} (q;q)_{\infty}\\
&&=(\delta_{\eta\gamma q^{s}} *_{\gamma} \delta_{\eta\gamma
q^{t}})(\theta\gamma q^{l}).\\
\end{eqnarray*}
Hence two discrete delta functions commute if and only if their supports have
the same signature. Therefore if 
two functions have support strictly contained in the same half line, we can
formally show that they commute by writing the two
functions as a sum of  discrete delta functions. For instance, if $f$ and $g$
are functions defined on $L(\gamma)^+:=\{q^{r}\gamma\,|\,r\in{\bf Z}\}$ then we
may write  formally $$f=\sum_{k=-\infty}^\infty f(q^{k}\gamma)\delta_{\gamma
q^{k}}\quad{\hbox{and}}\quad g=\sum_{l=-\infty}^\infty
g(q^{l}\gamma)\delta_{\gamma q^{l}}$$  so that their convolution product has
also support in $L(\gamma)^+$ and formally:
\begin{eqnarray*}
&&(f*_{\gamma}g)(\gamma q^{p})=\sum_{k=-\infty}^\infty\sum_{l=-\infty}^\infty
f(q^{k}\gamma)g(q^{l}\gamma)(\delta_{q^{k}\gamma}*_{\gamma}\delta_{q^{l}\gamma})(q^{p}\gamma)\\
&&=\sum_{k=-\infty}^\infty\sum_{l=-\infty}^\infty
f(q^{k}\gamma)g(q^{l}\gamma)(\delta_{q^{l}\gamma}*_{\gamma}\delta_{q^{k}\gamma})(q^{p}\gamma)
=(g*_{\gamma}f)(\gamma q^{p})\\
&&=(1-q)(q;q)_{\infty}(q^{p}\gamma)\sum_{h=0}^\infty\sum_{t=0}^\infty
{{(Q^pf)(q^{h}\gamma)}\over{(q;q)_{h}}}{{Q^pg)(q^{t}\gamma)}
\over{(q;q)_{l}}} q^{(h+t)+ht}.
\end{eqnarray*} Clearly $f$ and $g$ could not be analytic in a
neighbourood of zero unless they are zero on $q^k\gamma$ for $k\ge k_0$. On the
other hand, if $\eta\not=\epsilon$ the two products are different since the
support of the product will be contained in the half line having the same
signature  as the  discrete delta function on the right hand side of the
product.\\ 
\begin{remark}\rm $q$-distributions and their $q$-Fourier transform have been
studied in \cite{OR}. Olshanetski and Rogov define regular $q$-distributions as
those  distributions $D(\psi)$ for which there is a function $\psi$ such that 
$D(\psi)(f)=\int_\ga {\overline\psi}f$ for $\ga=1$. In particular,
$D(\delta_{\eta q^k\ga})(f)=(1-q) q^k\ga f(\eta q^k\ga)$. One can check in
this case that  $$D(\pa\delta_{\eta q^k\ga})(f)=-D(Q\delta_{\eta q^k \ga})(\pa
f)=-q^{-1}D(\delta_{\eta q^k\ga})(Q^{-1}\pa f)$$
Observe that the  ($q$-)regular distributions defined by
$(1-q)^{-1}q^{-k}\delta_{\eta q^k \ga}$ act as classical distributions on
test functions vith real values, and their limit for $k\to\iy$ is the ordinary
distribution given by Dirac's delta.\hfill$\spadesuit$\end{remark}

\section{The family $\MM_s$}\label{families}
For any $s>0$ let $\MM_s$ be the family  of
functions of the form $F=f\,e_{q^2}(-X^2)$ where
$f(x)=\sum_{l=0}^\infty a_lx^l$ with $|a_l|\le C s^l q^{\half{l^2}}$ for
some $C>0$.\\
It can be shown that if $s<q^{-\half}$, such an $f$ may also be written 
as  $f(x)=\sum_{l=0}^\infty c_l{\tilde h}_l(x;q)$ with $|c_l|\le C'
s^l q^{{l^2}\over2}$ for some $C'>0$ {\em and same $s$}. \\
Here, ${\tilde
h}_{l}(x;q) =(q;q)_k\sum_{l=0}^{[k/2]}{{(-1)^lq^{-2lk+2l^2+l}x^{k-2l}}
\over{(q^2;q^2)_l\,(q;q)_{k-2l}}}$
(see \cite{KoSw})
 are the discrete $q$-Hermite II polynomials.  Viceversa, if
$f(x)=\sum_{l=0}^\infty c_l{\tilde h}_l(x;q)$ with $|c_l|\le C' s^l
q^{{l^2}\over2}$ for some $C'>0$ and some $s\in(0,\,q^{-{1\over2}})$ then
$f(x)= \sum_{l=0}^\infty a_lx^l$ with $|a_l|\le C s^l q^{{l^2}\over2}$ for
some $C>0$. This is achieved by means of formulas (8.9) and (8.17) in
\cite{Koo}, and the estimate in the proof of Theorem 6.5 in
\cite{giotom}. Clearly $\MM_s\subset\HS$ for every $s>0$ and 
$\MM_s\subset\MM_r$ if $s<r$. \\ 
\begin{remark}\rm In \cite{ramis}, where $q>1$, power series of the type
above described are called $q$-Gevrey of order $-1$ and finite type.
$q$-Gevrey series were first introduced in \cite{bezi}, but only for
positive type. It is shown
in Proposition 2.1 in \cite{ramis} that the  above conditions on the
coefficients of a power series $f$ imply that $f$ has $q$-exponential growth 
of order $1$ and finite type. In particular, For $F=f\,
e_{q^2}(-X^2)\in\MM_s$,  $f$ will have  order $-1$ and finite type smaller or
equal to $sq^{\half}$. \hfill$\spadesuit$\end{remark} 
By a simple computation
one sees that   $\MM_s\subset\infigamma$  for every $s>0$ and for $\ga>0$.
Indeed, \begin{eqnarray*}&&\int_{\gamma}|f\,e_{q^2}(-X^2) X^e|\le
C\sum_{n=0}^\infty s^n q^{{n^2}\over2}\int_{\gamma}|X|^{n+e}e_{q^2}(-X^2)\\
&&\le C'\sum_{n=0}^\infty  q^{{{n^2}\over 2}-{{(n+e)^2}\over
2}+{{-n-e}\over2}+{{1\over4}(n+e)^2}}s^nb^{n+e}<\infty \end{eqnarray*} for
some positive constants $C,\,C'$ and $b$. Here is  used that $e_{q^2}(-X^2)$ is
of strict left type  $1/2$, as it was shown in \cite{giotom}. Hence, 
${\MM}_s\subset \HS\infigamma$ for every $\gamma\in(0,\,1)$.  \\ However,  the
elements of ${\MM}_s$ do not belong to $\II^\com_\ga$ in general. Indeed, let
$M(x)=e_{q^2}(-x^2)\sum_{n=0}^\infty q^{2n^2} s^nx^{2n}$ for some $s<1$.  One
computes $$\mu_{2e,\gamma}(M)=c_q(\gamma)q^e\sum_{n=0}^\infty
q^{(n-e)^2}s^n(q;q^2)_{n+e}$$ where $c_q(\ga)$ is as in formula
(\ref{eq-moment}). Since $$\sum_{n=0}^\infty
q^{(n-e)^2}s^{n-e+e}(q;q^2)_{n+e}\ge (q;q^2)_{\infty} s^e\sum_{n=e}^\infty
q^{(n-e)^2}s^{n-e}=(q;q^2)_{\infty} s^e\sum_{p=0}^\infty q^{p^2}s^{p}$$
$|\mu_{2e,\gamma}(M)|\ge Cq^es^e$ for some constant $C$ hence
$M\not\in\II^\com_\ga$.\\  If $g\in \II^\com-\ga$ and $F\in{\cal M}_s$
for $s\in(0,\,q^{-{1\over2}})$ it makes sense to compute $g*_{\gamma}F$, 
which will belong to  $\HS\II_{\ga}^\iy$  for every $\gamma<1$ by Proposition
4.4  and Lemma 3.4 in \cite{giotom}. \begin{proposition} For every
$s\in(0,\,q^{-\half})$, ${\cal M}_s$ is a left module for $\HD\II^{\cal
E}_\ga$,  its subalgebras  and their quotient   by their respective commutator
ideals. \end{proposition}   \Proof   By $(8.28)$ in \cite{Koo}
\begin{equation}\label{hermitepartial}\partial^t\bigl({\tilde
h}_l(x;q)e_{q^2}(-x^2)\bigr)={{(-1)^tq^{lt+{{t^2-t}\over2}}}\over
{(1-q)^t}}{\tilde h}_{l+t}(x;q)e_{q^2}(-x^2).\end{equation} Hence, for
$g\in\HD\II_\gamma^{\cal E}$ and $F(x)=\sum_{t=0}^\infty a_t {\tilde h}_t(x;q)
e_{q^2}(-x^2)$ one can compute  \begin{eqnarray} (g*_{\gamma}F)(x)
&=&\sum_{e=0}^\infty{{q^{{e^2+e}\over2}\int_{\gamma}g\,X^e}\over{(q;q)_e}}
\sum_{p=e}^\infty a_{p-e}q^{(p-e)e+{{e^2-e}\over2}}{\tilde
h}_{p}(x;q)e_{q^2}(-x^2)\\
&=&\sum_{p=0}^\infty\Biggl[\sum_{e=0}^p{{q^{ep}\int_\gamma
g\,X^e}\over{(q;q)_e}}a_{p-e}\Biggr]{\tilde
h}_p(x;q)e_{q^2}(-x^2)\label{propro} \end{eqnarray} where we could interchange
summations by dominated convergence, using the estimate 
\begin{equation}\label{giovanna} {q^{k^2-k\over 2}\over(q;q)_k}\,|\tilde
h_k(x;q)| \le{{(1-q)^{-1}\max(1,|x|)}\over{(q;q)_k}}\sum_{p=0}^\iy
{{q^{2p^2-p}|x|^{2p}}\over{(q^2;q^2)_p}}\end{equation} which was obtained in
the proof of Theorem 6.5 in \cite{giotom}. One has
$$\Biggl|\sum_{e=0}^p{{q^{ep}\int_\gamma
g(x)x^e}\over{(q;q)_e}}a_{p-e}\Biggr|\le
Bs^pq^{p^2\over2}\sum_{e=0}^\infty{{|\mu_{e,\ga}(g)|s^{-e}q^{-\half
e}}\over{(q;q)_e}}= {\tilde D} q^{{p^2}\over2} s^p$$ for nonnegative constants
$B$ and $\tilde D$. Hence $g*_\ga{\cal M}_s\subset{\cal M}_s$. The fact that
$(f*_\ga g)*_\ga F=f*_\ga(g*_\ga F)$ holds already for $g\in\HD$. Last
statement follows by equation (\ref{generating}).\hfill$\square$\smallbreak 

In particular, for every
$s\in(0,\,q^{-{1\over2}})$ the familes ${\cal M}_s$ are left modules for all
algebras of functions of left type, strict left type, etcetera. 
\begin{corollary}\label{ideal} For every $s\in(0,q^{-{1\over2}})$  and for
every $\gamma$ the space $\II^\com_\gamma\cap{\cal M}_s$ is an algebra and
$\II^\om_\gamma\cap{\cal M}_s$ a subalgebra. For $s,\,r\in(0,q^{-{1\over2}})$
with $s<r$, $\II^\com_\ga\cap{\cal M}_s$ is an ideal of  $\II^\com_\ga\cap{\cal
M}_r$ and a module over $\HS\II^\com_\ga\cap{\cal
M}_r/[\HS\II^\com_\ga\cap{\cal M}_r,\,\HS\II^\com_\ga\cap{\cal M}_r]_*$.
Analogous  statements hold when we replace everywhere  the upper index $\cal E$
by $\omega$.\hfill$\square$\end{corollary}  \smallbreak
It follows by the particular structure of the elements of $\MM_s$ that for a
fixed $s$, $f\in\MM_s\cap\II^\com_\ga$ if and only if
$f\in\MM_s\cap\II^\com_{\ga'}$ for $\ga,\ga'<1$ and that the two algebras are
isomorphic because $f*_\ga g={{c_q(\ga)}\over{c_q(\ga')}}f*_\ga' g$.
Therefore we could even remove the lower index $\ga$ from $\II_\ga$,
$*_\ga$, $\mu_\ga$ etcetera. 

\begin{remark}\label{closed}\rm By formula (\ref{hermitepartial})
it follows by direct computation that every ${\cal M}_s$ with
$s\in(0,\,q^{-{1\over2}})$ is closed under $q$-differentiation. \\
Indeed
if $F=(\sum_k c_k {\tilde h}_k(x;q))\,e_{q^2}(-X^2)$ with $|c_k|\le C s^k
q^{\half k^2}$, then \\
$\partial^r F=(\sum_kd_k  {\tilde
h}_k(x;q))\,e_{q^2}(-X^2)$ with  $|d_k|\le (C q^{-\half r}
s^{-r}(1-q)^{-r})s^k q^{\half k^2}$. \\
However, the family is not closed under
ordinary multiplication by $x$ or under $q$-shift $Q$. In general, there hold
only the weaker formulas, for $s<q^{\half}$: $X\MM_s\subset \MM_{sq^{-1}}$ and 
$Q\MM_s\subset \MM_{sq^{-1}}$ as it follows by direct computation, using
$e_{q^2}(-x^2)={1\over{(-x^2;q^2)_\iy}}$.\hfill$\spadesuit$ \end{remark} 
\begin{example}\rm  Let $p(x)=\sum_{n=0}^M a_nx^n\in\CC[x]$. Then
$p(X)\,e_{q^2}(-X^2)\in {\cal M}_s$ for every $s$ by taking
$C=\max_n(|a_n|)q^{-{{M^2}\over2}}s^{-M}$. Clearly
$p(X)\,e_{q^2}(-X^2)$ is of strict left type $1\over2$ 
since $e_{q^2}(-X^2)$ is (see Example 3.2 in \cite{giotom}).\hfill$\spadesuit$
\end{example} \begin{example}\label{giemme}\rm
For $m\in\Zplus$ let
\begin{eqnarray}
&&g_m(x):=
e_{q^2}(-x^2)\,{}_0\phi_1(-;q^{1+2m};q^2,-q^{1+2m}x^2)\nonumber\\
&&\qquad\quad=e_{q^2}(-x^2)
\sum_{r=0}^\iy{(-1)^r q^{2r^2-r} q^{2mr} x^{2r} \over
(q^{1+2m};q^2)_r (q^2;q^2)_r}\,.\label{functions}
\end{eqnarray}Those functions were constructed first in Example 6.7 in \cite{giotom}.
One sees immediately that $g_m\in\MM_{q^{m-\half}}$ for every $m\in\ZZ_{\ge1}$.
By the  formula after  (2.14) in \cite{KooSw},
$g_m(x)={}_2\phi_1(0,0;q^{2m+1}|q^2,-x^2)$ and the $g_m$'s are
related to Jackson's $q$-Bessel function ${\cal J}^{(1)}_{\al}(2x;q^2)$ by:
\begin{equation}\label{bessel}g_m(x)=
{{(q^2;q^2)_\iy}\over{(q^{2m+1};q^2)_\iy}}x^{-m+\half}{\cal
J}^{(1)}_{m-\half}(2x;q^2)\end{equation} (see \cite{KoSw}, \cite{KooSw} and
references therein). It was shown in \cite{giotom} that  for every
$m\in\Zplus$, $g_m\in\HS\II^\om_{\ga,\al}$ for every $\alpha>0$
because the $g_m$'s are even and $\mu_{2k}(g_m)=0$ for every $k\ge m$. 
The $q$-moments were explicitely
computed: the odd ones are always zero and 
$$\mu_{2k,\ga}(g_m)=c_{q}(\ga)q^{k^2+k}(q;q^2)_k\,{{(q^{2m-2k};q^2)_{\iy}}
\over{(q^{1+2m};q^2)_{\iy}}}$$ where $c_q(\ga)$ is as usual. Hence 
 convolution from the left by some
$g_m$ is  equivalent to a genuine $q$-differential operator of order $2m$.\\
It was shown in \cite{giotom} that $g_0$ and $g_1$ are not of strict
left type. Using the three term recurrence relation for Jackson's
$q$-Bessel function in Exercise 1.25 in \cite{gasper}, one sees that
the same holds for all $g_m$'s. In particular one could use the same
three term recurrence relation in order to show that the $g_m$'s are
of left type for all $\al$ once this is shown for $g_0$ and $g_1$.\\
Clearly $u_\ga:={{g_1}\over{\int_\gamma g_1}}$ is a left unit for
$\II^\com_\gamma\cap{\cal M}_s$ for every $s\in[q^\half,\,q^{-\half})$, hence
for  $\HS\II^\om_\gamma\cap{\cal M}_s$ for $s$ ranging in the same set.  Later
we will see that $u_\ga$ is also a right unit too, for $s\in[q^\half,\,1)$.
Moreover,  since by Lemma 4.3 in \cite{giotom}
\begin{equation}\label{unit}(\pa^k u_\ga)*_\ga f=\partial^k(u_\ga*_\ga
f)=\partial^k f\quad\quad\forall k\in\Zplus\end{equation} 
convolution by  $\pa^k u_\ga$  from the left coincides with applying
$\pa^k$.\hfill$\spadesuit$ \end{example}

For functions in $\MM_{s}$ with $s<q^{-\half}$ it is possible to
improve Lemma \ref{domain}.\\ 
 \begin{proposition}\label{productzero} Let $g\in
\II^\com_\ga$ and let $F\in\MM_s$ with $s<{q^{-\half}}$. Then $g*_\ga F=0$ iff 
$\mu_{\ga}(g)=0$  and/or $F=0$.  In particular, if $g$
belongs to the annihilator of some nonzero function $F\in\MM_s$
with $s<{q^{-\half}}$, then  $g*_\ga f=0$ for every $f$ for which the product
is defined. In particular the representation of
$\HD\II_\ga^\com/[\HD\II_\ga^\com,\, \HD\II_\ga^\com]_*$, of
$\HD\II_\ga^\om/[\HD\II_\ga^\om,\, \HD\II_\ga^\om]_*$, and of
$\HS\II^\som_{\ga,>\half}$ on $\MM_s$ are faithful. 
\end{proposition} \Proof  Use formula (\ref{propro})  and independence
of $q$-Hermite polynomials in order to show that $g*_\ga F=0$ if and only if
either $\mu_\ga(g)=0$ of $F=0$. The last statement
follows by Lemma 6.1 in \cite{giotom}.\qed    
\begin{remark} \rm Observe  that
the result in Proposition \ref{productzero} does not necessary hold in
general. Namely,  the function with zero moments for a product equal to
zero might not be the left one. Take for instance
$e_{q^2}(-q^2X^2)*_\ga e_q(iX)=0$ that was computed in Remark \ref{domicommu}.
\hfill$\spadesuit$\end{remark}

We will investigate
commutativity. The following Lemma was communicated to me by T. Koornwinder.
\begin{lemma}\label{tom} Let $s\in(0,\,1)$ and let $F=f\,e_{q^2}(-X^2)\in
{\cal M}_s$ and even. Suppose also that for some $\gamma>0$ we have: 
$\int_{\gamma}X^{2k}\,F=0$ for all $k\in{\bf Z}_{\ge0}$. Then $F=0$.\end{lemma}
{\bf Proof}: Let $F=f\,e_{q^2}(-X^2)\in {\cal M}_s$ 
with  $f (x)=\sum_{n}a_{2n}x^{2n}$ and 
$|a_{2n}|\le C\,s^n\,q^{2n^2}$ for some $C>0$ and some $s\in(0,1)$.
It  will be justified  by dominated convergence that
$$
\int_{\gamma}|f(x)|^2\,e_{q^2}(-x^2)=
\lim_{m\to\infty}\int_\gamma
\left(\sum_{k=0}^m \overline{a_{2k}}\,x^{2k}\right) f(x)\,
e_{q^2}(-x^2)\,d_qx=0.
$$
Hence $f(\gamma q^k)=0$ for all $k\in{\bf Z}$. So $f\in{\cal E}$ vanishes on a set with limit point 0. 
Hence $f=0$ identically.
For the proof of the dominated convergence note that, for $x>0$,
$$
\left(\sum_{n=0}^\infty|a_{2n}|\,x^{2n} \right)^2\le
C\sum_{n=0}^\infty x^{2n}s^n
\left(\sum_{k=0}^n q^{2k^2}q^{2(n-k)^2}\right)
$$
$$
\le
C'\sum_{n=0}^\infty x^{2n}s^n q^{n^2}\le C''E_{q^2}(qsx^2).
$$ for constants $C$, $C'$ and $C''$.
Hence
$$
\sum_{k=-\infty}^\infty|f(\gamma q^k)|^2\,e_{q^2}(-q^{2k}\gamma^2)\,q^k\le
D\sum_{k=-\infty}^\infty
{(-\gamma^2;q^2)_k\,q^k\over
(-qs\gamma^2;q^2)_k}<\infty.
$$ where $D$ is a nonnegative constant.
This completes the proof.\hfill$\square$\smallbreak One  can extend the above
Lemma to prove an analogue result for all functions in ${\cal M}_s$ for $s<1$.
\begin{lemma}\label{tom-odd} Let 
$s\in(0,\,1)$ and let $F=f\,e_{q^2}(-X^2)=\sum_{n=0}^\infty
a_nx^ne_{q^2}(-x^2)\in{\cal M}_s$. If  for some $\gamma>0$ one has
$\int_{\gamma}F\,X^k=0$ for every $k\in{\bf Z}_{\ge0}$ then $f=F=0$.
\end{lemma}
{\bf Proof: }We write $f=f_0+f_1$ where $f_0$ (resp. $f_1$) is the even
(resp. odd) part  of $f$. Then,
$\int_\gamma F\,X^{2k}=\int_\gamma f_0\,e_{q^2}(-X^2) X^{2k}=0$ for every
$k\ge0$ so, by Lemma \ref{tom} $f=f_1$. Then
$\partial f\in\MM_s$ by Remark \ref{closed} and it is
even. $\int_\ga \pa f =0$ and $\int_\ga X^k\pa
f=-q^{-k}[k]_q\int_\ga X^{k-1}f=0$ for every $k\ge1$.
Apply Lemma \ref{tom} to get the statement.\hfill$\square$\smallbreak
In particular we have:
\begin{corollary} For every $\gamma>0$ and for $s\in(0,\,1)$,
$\II^\com_\gamma\cap{\cal M}_s$  is a commutative algebra under the convolution
product.  If moreover, $s\in[q^\half,\,1)$,  $\II^\com_\gamma\cap{\cal M}_s$
is unital. \end{corollary} \Proof Since by formula (\ref{generating}) 
$\mu_{\ga}(f*_\ga g-g*_\ga f)=0$, the first statement follows by  Lemma
 \ref{tom-odd}. The second statement follows by
commutativity, formula (\ref{unit}) for $k=0$ and the fact that 
$u_\ga\in\II^\om_\ga\cap\MM_{q^\half}$.
\hfill$\square$\smallbreak It follows that the subalgebras
$\II^\om_\gamma\cap{\cal M}_s$, $\II^\somE_\gamma\cap{\cal M}_s$ and
$\II^\som_\gamma\cap{\cal M}_s$  are also commutative for $s\in(0,\,1)$, and
that $\II^\om_\gamma\cap{\cal M}_s$  is unital for $s\in[q^{\half},\,1)$.
\begin{corollary} For every $\gamma>0$ and for $s\in[q^\half,\,1)$,
$\II^\om_\gamma\cap{\cal M}_s=\II^\om_\ga\cap\MM_{q^\half}$, and
$\II^\com_\gamma\cap{\cal M}_s=\II^\com_\ga\cap\MM_{q^\half}$.\end{corollary}
\Proof  By commutativity and Corollary \ref{ideal},
$\II^\om_\ga\cap\MM_{q^\half}$  and $\II^\com_\ga\cap\MM_{q^\half}$  are
bilateral ideals of $\II^\om_\ga\cap\MM_{q^\half}$ and
$\II^\com_\ga\cap\MM_{q^\half}$ respectively, and they contain $u_\ga$. Hence
they coincide. \hfill$\square$\smallbreak   \begin{example}\rm  Functions of
the form $p(X)\,e_{q^2}(-X^2)$ where $p(X)$ is a polynomial function, form a
commuting family of functions. They all belong to $\II_\ga^\som$
for every $0<\ga<1$. However this class of functions is
not closed under convolution product.  Commutativity can also be checked
directly as follows. Take  ${\tilde h}_l(x;q)e_{q^2}(-x^2)$  for $l\in\Zplus$
as basis for the above space.   By formulas $(8.17)$  and $(8.14)$ in
\cite{Koo} $$\int_{\gamma}e_{q^2}(-x^2)x^p{\tilde h}_r(x;q)=
\cases{c_q(\gamma){{(q;q)_{r+2k}q^{-2rk-k^2-r^2}}\over{(q^2;q^2)_k}}& if
$p-r=2k\ge0$;\cr 0& otherwise.\cr}$$  Hence   \begin{eqnarray*}&&({\tilde
h}_r(X;q)e_{q^2}(-X^2)*_{\gamma}{\tilde
h}_l(X;q)e_{q^2}(-X^2))(x)\\
&&=c_q(\gamma)q^{lr}\sum_{k=0}^\infty
{{q^{3k^2}q^{2(r+l)k}}\over{(q^2;q^2)_k}}{\tilde
h}_{l+r+2k}(x;q)e_{q^2}(-x^2)
\end{eqnarray*}  Since the last expression is
symmetric in $l$ and $r$ commutativity holds, but the product will no
longer be  a
polynomial times $e_{q^2}(-X^2)$.\hfill$\spadesuit$
\end{example} 
\begin{example} \rm The functions $g_m$ for $m>0$ defined in Example
\ref{giemme}  are a family of commuting functions in ${\cal
M}_{q^\half}$ by Lemma \ref{tom}. On the contrary, $g_0\not\in{\cal
M}_s$ for any $s<q^{-\half}$ (its coefficients grow exactly like
$q^{{n^2}\over2}q^{-\half n}$ as $n\to\iy$) and as it was shown in Example 6.7
in \cite{giotom} it does not commute with $g_1$, hence with $u_\ga$. On the
other hand one can conclude that
$$g_m*_{\gamma}g_n=\sum_{r=0}^{m-1}{{\mu_{2r,\ga}(g_m)}\over{(q;q)_{2r}}}\partial^{2r}g_n=g_n*_{\gamma}g_m=
\sum_{r=0}^{n-1}{{\mu_{2r,\ga}(g_n)}\over{(q;q)_{2r}}}\partial^{2r}g_m$$
for $m,\,n$ in ${\bf Z}_{\ge1}$.
Using
determinacy
of the  $q$-moment
problem, one has for $m\in{\bf Z}_{\ge1}$
$$g_m(x)={{(q^{2m};q^2)_\iy}\over{(q^{1+2m};q^2)_\iy}}\Biggl(\sum_{k=0}^\iy{{(-1)^k q^{2mk} q^{2k^2-k} 
{\tilde h}_{2k}(x;q)}\over{(q^2;q^2)_k}}\Biggr) e_{q^2}(-X^2)$$ since both
functions are in $\MM_{q^{m-\half}}\cap\II^\om_\ga$ and have the same
$q$-moments. Hence we have the equality:
$${}_0\phi_1(-;q^{1+2m};q^2,-q^{2m+1}x^2)={{(q^{2m};q^2)_\iy}\over{(q^{1+2m};q^2)_\iy}}  
\Biggl(\sum_{k=0}^\iy{{(-1)^k q^{2mk} q^{2k^2-k} {\tilde
h}_{2k}(x;q)}\over{(q^2;q^2)_k}}\Biggr)$$ for $m\ge1$.
\hfill$\spadesuit$\end{example} \begin{corollary}\label{annihilator} For every
$s\in(0,\,q^\half]$ and $r<q^{-\half}$ the representation of
$\II^\com_\gamma\cap{\cal M}_s$  on $\MM_r$ is
faithful. This implies, taking  $r=s$ that $\II^\com_\gamma\cap{\cal M}_s$ and its subalgebras  are
algebraic domains.\end{corollary} \Proof  By Lemma
\ref{productzero}  the annihilator of $\MM_r$ in $\II^\com_\gamma\cap{\cal
M}_s$  must be zero. \hfill$\square$\smallbreak 
\begin{corollary}\label{nucleo}
$\HD\II^\om_\ga/[\HD\II^\om_\ga,\,\HD\II^\om_\ga]_*$ is a unital algebra. The
same holds if we replace the upper indices $D$ by $S$ and/or $\omega$ by
$\com$.  Moreover for the kernel of the formal $q$-Fourier transform $\tilde
F_\ga$ one has: \begin{eqnarray*}&&Ker_{\HD\II^\om_\ga}(\tilde\FF_\ga)=\{
f\in\HD\II^\om_\ga\,|\, \mu_\ga(f)=0\}=[\HD\II^\om_\ga,\,\HD\II^\om_\ga]_*\\ 
&&Ker_{\HS\II^\om_\ga}(\tilde\FF_\ga)=\{ f\in\HS\II^\om_\ga\,|\,
\mu_\ga(f)=0\}=[\HS\II^\om_\ga,\,\HS\II^\om_\ga]_*
\end{eqnarray*} and the same if we replace everywhere $\omega$
by $\com$. \end{corollary} \Proof The proof will be for
$\HS\II^\om_\ga$, the other cases follow similarly. The
function $u_\ga$ is a left unit in $\HS\II^\om_\ga$ hence its projection on 
the commutative algebra $\HS\II^\om_\ga/[\HS\II^\om_\ga,\,\HS\II^\om_\ga]_* $
is a unit therein.  The first equality of the formula for the kernel in
$\HS\II^\om_\ga$ is clear.  Inclusion $\supseteq$ in the second equality 
follows by equation (\ref{generating}). Let $f\in\HS\II^{\om}_\ga$  be such
that $\mu_\ga(f)=0$.  Then $f=[u_\ga,\,f]_*$, hence
the other inclusion. \hfill$\square$ \smallbreak 
\noindent Observe that
$u_\ga$ is not a left unit on $\HD\II^\om_\ga$ since $e_q(iX)*_\ga
u_\ga=0$. \smallbreak 
\section{The functions $u_\ga$ and $G_{k,\ga}$ and
topology}\label{Gicappa}

In this section we  study of the functions $g_m$, defined in formula
(\ref{functions}) and their $q$-derivatives, as they are of particular interest
and  useful in order to prove plenty of results. \\ 
Let $k,\,m\in\ZZ_{\ge1}$ and let
$g_{k,m}:=\partial^{k}g_m$. By Remark \ref{closed},  $g_{k,m}\in
\II^\om_{\ga,\al}\cap\MM_{q^{m-\half}}$ for every $\al$ and for every
$\ga\in(0,\,1)$. Moreover, by  Lemma 4.1 in \cite{giotom}, we have
$$\mu_{l,\ga}(g_{k,m}) =\cases{0 & if $k+l$ odd,\cr 0 & if $l<k$,\cr
(-1)^k{{[l]_q!}\over{[k]_q!}}\mu_{l-k}(g_m)& otherwise.\cr}$$ In particular
for  $k=2p$ we have $\mu_{2j+1,\ga}(\pa^{2p}g_m)=0$ and
$\mu_{2j,\ga}(\pa^{2p}g_m)=0$ for every $j\ge m+p$. One can see that
$g_{2p,m}$ is not a multiple of $g_{m+p}$ for $m\ge1$ since
$\mu_{2j,\ga}(g_{m+2p})\not=0$ for  $j<p$.  On the other hand,
$\pa^{2k}g_0={{(-1)^k}\over{(1-q)^{2k}}}g_0$. This is checked using the fact
that $g_0(x)=\half(e_q(ix)+e_q(-ix))=cos_q(x)$ (see
Example 6.7 in \cite{giotom}). \\
Let us introduce the family of functions in
$\MM_{q^\half}\cap\II^\om_\ga$:
\begin{equation}\label{basis}G_{k,\ga}:={{(-1)^{k}g_{k,1}}\over{[k]_q!\int_\ga
g_1}}={{(-1)^{k}\pa^{k}u_\ga}\over{[k]_q!}}=
{{(-1)^{k}(q^2;q^2)_\iy\pa^{k} \biggl(x^{-\half}{\cal
J}_\half^{(1)}(2x;q^2)\biggr)}\over{(q^3;q^2)_\iy [k]_q!}}\end{equation}  By
direct computation one obtains: 
\begin{equation}\label{leggi}G_{l,\,\ga}(x)={{e_{q^2}(-x^2)}\over{c_q(\ga)(q;q)_l}}
\sum_{k=0}^\iy{{(-1)^k q^{\half(2k+l)^2+\half(2k-l)}{\tilde
h}_{2k+l}(x;q)}\over{(q^2;q^2)_k}}\end{equation} and $G_{l,\ga}$ differs from
$G_{l,\ga'}$ only by a multiplicity constant, since this is true for
$u_\ga$. Observe that
$g_1(x)x(1-q)^{-1}=\sin_q(x)={1\over{2i}}(e_q(ix)-e_q(-ix))$ so that in
particular the limit for $q\to1^{-}$ of $g_1((1-q)x) x$ is $\sin(x)$.\\ The
$G_{k,\ga}$'s belong to $\MM_{q^\half}\cap\II_{\ga,\al}^\om$ for every $\al$
and every $\ga$ since  $\mu_{r,\ga}(G_{k,\ga})=\delta_{r,k}$.
$G_{k,\ga}*_\ga f=(-1)^{k}{{\pa^{k}}\over{[k]_q!}}\pa^k f$ for every $f$ for
which the $q$-derivatives are defined, and \begin{equation}G_{k,\ga}*_\ga
G_{l,\ga}=\Bigl[{{k+l}\atop l}\Bigr]_qG_{k+l,\ga}.\end{equation}
The  $g_m$'s can be written  as linear combinations of the $G_{r,\ga}$'s using
determinacy of the $q$-moment problem in $\MM_{q^\half}$. Indeed 
$$g_m= \sum_{r=0}^m\mu_{2r}(g_m)
\,G_{2r,\ga}={{c_q(\ga)(q^{2m};q^2)_\iy}\over{(q^{1+2m};q^2)_\iy}}\sum_{r=0}^m
(-1)^rq^{2mr}(q;q^2)_r(q^{2-2m};q^2)_rG_{2r,\ga}.$$ 
since  both functions belong to
$\MM_{q^\half}\cap\II^\om_\ga$   and  have the same
$q$-moments. This provides another way to express Jackson's $q$-Bessel
functions using formula (\ref{bessel}). \\
The expansion of functions in
$\MM_{q^\half}\cap\II^\om_\ga$  in terms of the $G_{k,\ga}$'s can
be seen as  an approximation of functions  with respect to  a
suitable topology, i.e. the one determined by the multiplicity of a zero at
$x=0$ of the $q$-moment series $\mu_\ga$. This is the subject of tha
last part of this section. \\ 
Let $\ga>0$ be fixed. For any $f\in\II_\ga^\iy$ let
$m(f)$ be the miminum nonnegative integer for which $\mu_{r,\ga}(f)\not=0$
(i.e. the multiplicity of a zero at $t=0$ of $\mu_\ga(f)$ or the multiplicity
of a zero at $y=0$ of $\tilde\FF_\ga(f)$).\\ Then  put, for
$f,\,g\in\II_\gamma^\iy$, 
\begin{equation}\d(f,\,g):=e^{-\m(f-g)}=\d(g,\,f)\end{equation} so that for
every $f$ and $g$ there holds $0\le\d(f,\,g)\le1$. Let
$f,\,g,\,h\in\II_\ga^\com$. If  $\mu_\ga(f-h)=t^{\m(f-h)}\phi(t)$ and
$\mu_\ga(g-h)=t^{\m(g-h)}\psi(t)$ where $\phi$ and $\psi\in{\cal E}$ are such
that  $\psi(0)\phi(0)\not=0$, then 
$\mu_\ga(f-g)=t^{\min(\m(f-h),\,\m(g-h))}F(t)$ with $F\in{\cal E}$ and one
has $\m(f-g)\ge\min(\m(f-h),\,\m(g-h))$. Hence  \begin{equation}\d(f,\,g)\le
e^{-\min(\m(f-h),\,\m(g-h))}=\max\bigl(e^{-\m(f-h)},\,e^{-\m(g-h)}\bigr)\le
\d(f,\,h)+\d(h,\,g)\end{equation} Therefore $\d$ defines a metric  on those
subspaces of  $\II_\gamma^\com$ on which the $q$-moment problem is determined.
For instance, $\d$ is a metric on the spaces $\MM_s\cap\II^\com_\ga$ with
$s\in(0,\, q^\half]$, the algebras $\HS\II_{\ga,>\half}^\som$,
$\HD\II_\ga^\om/[ \HD\II_\ga^\om,\,\HD\II_\ga^\om]_*$ and
$\HD\II^\som_{\ga}/{\rm Ann}(u_\ga)$, where ${\rm Ann}$ denotes the
annihilator ideal.\\ Let $f,\,g$ and $h$ belong to  an algebra with respect
to the $q$-convolution product. Then by formula (\ref{generating}) one has 
 $\m(f*_\ga
g)=\m(f)+\m(g)$ so that  $\d(f*_\ga g,\,h*_\ga g)=\d(f,\,h)\d(g,\,0)\le
\d(f,\,h)$. Hence the $q$-convolution product  with $g$ for a fixed $g$ is
continuous with respect to this topology. In particular, convolution by
$\pa^ku_\ga$, i.e. $q$-differentation, is continuous, as can be seen by the
fact that $\m(\pa^k f)=\m(f)+k$.\\ 
We have:
\begin{lemma}\label{appa} Let $g\in\II^\iy_\ga$ be such that $\mu_{\ga}(g)$
converges absolutely at least on the closed disk centered at zero and with
radius $(1-q)^{-1}$.  Then the sequence of functions
$F_n(g):=\sum_{k=0}^n\mu_{k,\ga}(g)G_{k,\ga}$ converges to a well-defined
function $F$. If the radius of convergence of $\mu_\ga(g)$ is strictly greater
than $(1-q)^{-1}q^{-1}$ or if $\gamma<1$ and $\rho<(1-q)^{-1}2\gamma^{-1}$
then  $F\in\II^\iy_\ga$ and  $\mu_\ga(F)=\mu_\ga(g)$. \end{lemma} \Proof  
Define  $F_n:=\sum_{k=0}^n\mu_{k,\ga}(g)G_{k,\ga}$ for every $n\in\Zplus$.  By
formula (\ref{leggi})  and  the estimate (\ref{giovanna})  we have
\begin{eqnarray*}&&\Biggl|\mu_{k,\ga}(g)\,G_{k,\ga}(x)\Biggr|\\
&&\le
C\,\max(1,\,|x|)\,\bigl|e_{q^2}(-x^2)\bigr|
\,\bigl|\mu_{k,\ga}(g)\bigr|\,\Biggr[\sum_{p=0}^\iy{{q^{2p^2-p}|x^{2p}|}\over{(q^2;q^2)_p}}\Biggr]\, 
\Biggl(\sum_{l=0}^\iy {{q^{2l}}\over{(q^2;q^2)_l}}\Biggr)
\end{eqnarray*}
for some constant $C$.
Hence $F_n$ tends to a well-defined function $F$ as $n\to\iy$ if 
$\sum_{k=0}^\iy |\mu_{k,\ga}(g)|<\iy$. \\ Next we want to  prove that  if the
radius of convergence $\rho$ of $\mu_\ga(g)$ is stricly greater than
$(1-q)^{-1}q^{-1}$ or if $\gamma<1$ and $\rho>2(1-q)^{-1}\gamma^{-1}$ then
$F\in\II^\iy_\ga$.  If  $\rho>(1-q)^{-1}q^{-1}$ one uses the fact that
$u_\ga=\sum_{k=0}^\iy c_kx^ke_{q^2}(-x^2)$ with $|c_k|\le C_0 q^{\half
(k^2+k)}$ for some constant $C_0$ together with Remark \ref{closed} in order
to conclude that  for every $l\in\Zplus$ 
$$\int_\ga |G_{p,\ga}| \,|X^l|\le C_0
{{q^{-p}}\over{(q;q)_p}}\sum_{k=0}^\iy
q^{\half(k^2+k)}q^{-\half(l+k)^2-\half(l+k)}\nu_{l+k,\ga}(e_{q^2}(-X^2))$$
where $\nu_{k,\ga}$ is defined in formula (\ref{momenteq}). Since it can be
shown that $\nu_{r,\ga}(e_{q^2}(-X^2))\le 2q^{-{1\over4}}c_q(\ga)
q^{{r^2}\over4}$ one obtains: $$\int_\ga|F X^l|\le
C_1q^{-\half(l^2+l)}\sum_{p=0}^\iy
{{q^{-p}|\mu_{p,\ga}(g)|}\over{(q;q)_p}}\sum_{k=0}^\iy
q^{{(l-k)^2}\over4}<\iy$$ hence $F\in\II^\iy_\ga$ and it is clear then that
$\mu_\ga(F)=\mu_\ga(g)$. If $\rho>2(1-q)^{-1}\gamma^{-1}$ and $\gamma<1$ we
may use Lemma 3.5 in \cite{giotom} in order to show that for every $l\in\Zplus$
$$\int_\ga |G_{p,\ga}|\,|X^l|\le
{1\over{[p]_q!}}\biggl(\int_\ga|X^l||u_\ga|+r^l\,
B\biggr){{2^p}\over{\gamma^p(1-q)^p}}$$ for some constants $B>0$ and
$r\in(\gamma,\,1)$. Then the proof follows as in the previous
case.\hfill$\square$\smallbreak

Observe that it follows by the proof of the
above Lemma that if  $g$ is any function of $\II^\iy_\ga$ for which
$\mu_{\ga}(g)$ has a big enough radius of convergence, then $g*_\ga u_\ga$ is
well defined. \begin{lemma} \label{construct2} Let $g\in\II^\iy_\ga$ be such
that $\mu_{\ga}(g)$ converges absolutely at least on the closed disk centered
at zero and with radius $(1-q)^{-1}$.  Then the function
$F:=\sum_{k=0}^\iy\mu_{k,\ga}(g)G_{k,\ga}$ belongs to $\MM_s$ for some
$s<q^{-\half}$. If the radius of convergence $\rho$ of $\mu_\ga(g)$ is strictly
greater than $(1-q)^{-1}q^{-1}$ then  $F\in\MM_{q^{\half}}$.  \end{lemma} 
\Proof  By dominated convergence 
$$F(x)={{e_{q^2}(-x^2)}\over{c_q(\ga)}}\sum_{p=0}^\iy q^{\half(p^2-p)}\tilde
h_p(x;q)\sum_{ k=0}^{[{p\atop
2}]}{{\mu_{p-2k,\ga}(g)q^{2k}}\over{(q;q)_{p-2k}(q^2;q^2)_k}}$$ 
Therefore the coefficients $c_p$ of the expansion of  $F\, E_{q^2}(X^2)$ with
respect to the discrete $q$-Hermite II polynomials are majorized by
$q^{\half(p^2-p)}\sum_{k=0}^{[{p\atop
2}]}|\mu_{p-2k,\ga}(g)|q^{2k}$ times some constant. If $\rho>(1-q)^{-1}$ then
$|\mu_{l,\ga}(g)|=O(a^l)$ for $l\to\iy$  for some $a\in(0,\,1)$. One
can always assume that $a\in(q,\,1)$. Then $|c_p|\le C q^{\half p^2}
(q^{-\half }a)^p$ for some constant $C$.  If $\rho>(1-q)^{-1}q^{-1}$ then $a$
can be chosen in $(0,\,q)$. In that case $|c_p|\le C'
q^{\half (p^2-p)} q^p$ for some constant $C'$.  \hfill$\square$\smallbreak

 \begin{corollary} Let $g\in\MM_s$ with $s<1$ be such that
$\mu_{\ga}(g)$ converges absolutely at least on the closed disk centered at
zero and with radius $(1-q)^{-1}q^{-1}$.  Then $g$ can be approximeted by
finite linear combinations of the $G_{k,\ga}$'s.\end{corollary} \Proof  By the
above results and Lemma \ref{tom-odd} there follows that
$F=g$.\hfill$\square$\smallbreak We have just seen that we can approximate
various classes of  functions by means of $u_\ga$ and its $q$-derivatives, as
in classical distribution theory one approximates generalized functions by the
delta functions and its derivatives. \\  Moreover, the methods used in the
proof of Lemma \ref{appa} and Lemma \ref{construct2} show that, if
$f\in\II^\iy_\ga$ is such that $\mu_\ga(f)$ has a good behaviour, there exists
a function $F\in\II^\com_\ga\cap\MM_{q^\half}$ such that
$\mu_\ga(f)=\mu_\ga(F)$. This is quite an interesting result because $F$
belongs to a space that  does not depend essentially on $\ga$. Moreover,
Lemma's \ref{appa} and \ref{construct2} provide a  {\em constructive} way to
associate to a n entire $q$-moment series a unique function in
$\MM_{q^\half}\cap\II^\com_\ga$. This proves the following theorem
\begin{theorem}\label{isomorfi} The elements of 
$\MM_{q^\half}\cap\II_\ga^\com$ form a set of representatives of the quotients
 $\HD\II_\ga^\com/[ \HD\II_\ga^\com,\,\HD\II_\ga^\com]_*$ and 
$\HS\II_\ga^\com/[ \HS\II_\ga^\com,\,\HS\II_\ga^\com]_*$. Those algebras
are all isomorphic. The same result holds if we replace everywhere the upper
index $\cal E$ by $\om$. The projection modulo the commutator ideal is given
in all cases by $f\mapsto f*_\ga u_\ga$.\end{theorem}  \Proof The bijection is
clear by the discussion above. The fact that it is an algebra isomorphism
follows from the fact that in those algebras the product is determined by
$\mu_\ga$.\hfill$\square$\smallbreak \noindent Observe that  even though
$\HS\II_\ga^\om$ is strictly contained in $\HD\II^\om_\ga$ and
$[\HS\II_\ga^\om,\,\HS\II_\ga^\om]_*$ is strictly contained in
$[\HD\II_\ga^\om,\,\HD\II_\ga^\om]_*$ (the function $\delta_{q^{-1}\gamma}$
belongs to $\HD\II^\om_\ga$ but not to $\HS\II^\om_\ga$), Theorem
\ref{isomorfi} states that for every function $f$ in $\HD\II_\ga^\om$ there is
always a function $f'\in \HS\II_\ga^\om$ such that
$\mu_\ga(f-f')\equiv0$.\smallbreak Observe also that the above result
together with Corollary \ref{nucleo} imply that the kernel of $\tilde
F_\ga$ on $\HD\II_\ga^\com$ does not depend on $\ga$ essentially.
\begin{remark}\rm On $\II_\ga^\com\cap\MM_{q^\half}$ and
$\II_\ga^\om\cap\MM_{q^\half}$   one can define the operators $\tilde X$ and
$\tilde Q$ as: ${\tilde X}.f=(Xf)*_\ga u_\ga$ and ${\tilde Q}.f=(Qf)*_\ga
u_\ga$.  It follows by formula (\ref{leibniz}) that $\tilde X$ acts as a
$q$-derivation. In particular, ${\tilde X}.G_{k,\ga}=q^{-k}G_{k-1,\ga}$  for
$k\ge1$ and ${\tilde X}.u_\ga=0$ as for the classical delta function (the unit
with respect to the convolution). As for the $q$-shift operator: ${\tilde
Q}.G_{k,\ga}=q^{-k-1}G_{k,\ga}$. Hence ${\tilde Q}{\tilde X}=q{\tilde
X}{\tilde Q}$.\end{remark}  Lemma's \ref{appa} and
\ref{construct2} can be generalised. \begin{lemma} \label{construct3} Let $h\in\MM_s$ for $s<
q^{-\half}$ and let $g\in\II^\iy_\ga$ be such that $\mu_{\ga}(g)$ converges
absolutely at least on the closed disk centered at zero and with radius
$(1-q)^{-1}q^{-\half}s^{-1}$.  Then $g*_\ga h$ is a well-defined function in
$\MM_s$.\end{lemma} \Proof  This Lemma generalises the result of Lemma
\ref{construct2} where $h=u_\ga$. One proves it similarly writing $h$ as
$e_{q^2}(-X^2)\,f$, expanding $f$ with respect to the discrete $q$-Hermite II
polynomials, and using the majorization of the coefficients of the expansion
of $\pa^k h$ given in Remark \ref{closed}.\hfill$\square$\smallbreak Let
$s<q^{-\half}$. By $\II_\ga^{\rho,s}$ we denote the space of functions in
$\II^\iy_\ga$ for which $\mu_\ga(f)\in\HD$ and has a radius of convergence
greater than $(1-q)^{-1}q^{-\half}s^{-1}$.\\
\begin{corollary}\label{new-algebra} For $s< q^{-\half}$, 
$\MM_s\cap\II_\ga^{\rho,s}$ is an algebra and equation (\ref{generating})
holds for functions in $\MM_s\cap\II_\ga^{\rho,s}$. If $s<1$ the algebra is
commutative. If $s=q^{\half}$, $\MM_{q^\half}\cap\II_\ga^{\rho,{q^\half}}$ is 
unital and all  $\MM_s\cap\II_\ga^{\rho,s}$ coincide for $s\in[q^\half,\,1)$.
\end{corollary} \Proof  By Lemma \ref{construct3} the product  of two
functions  in $\MM_s\cap\II_\ga^{\rho,s}$ is well-defined. One shows by
dominated convergence that associativity holds and one shows similarly to the
proof of Lemma \ref{appa} that  $\mu_\ga(f*_\ga g)$ is well-defined and that
it is equal to $\mu_\ga(f)\mu_\ga(g)$. By Lemma \ref{tom-odd} 
$\MM_s\cap\II_\ga^{\rho,s}$ is commutative, if $s<1$ and
$u_\ga\in\MM_{q^{\half}}\cap\II^{\rho,
q^{\half}}_\ga$.\hfill$\square$\smallbreak

\section{Convolution and Fourier transform}\label{convofourier}

In
this section I  apply the results of Section \ref{Gicappa} in order to extend
Koornwinder's inversion results for the $q$-Fourier transform to be found in
\cite{Koo}. At the end of the Section I  shall also prove  anaytically  a
result on the relation between $q$-convolution and $q$-Fourier transform that 
was proved in  \cite{KeMa} in a different context (braided and with bosonic
integral). \\
As it is stated in Section 2 the formal $q$-Fourier transform $\tilde F_\ga$
is defined on $\HD\II^\com_\ga$ and the $q$-Fourier transform is defined on
$\HD\II^{s\com}_\ga$. We have seen in Corollary \ref{nucleo} and Theorem
\ref{isomorfi} that 
$$\MM_{q^\half}\cap\II^\om_\ga\simeq\HS\II^\om_\ga/[\HS\II^\om_\ga,\,\HS\II^\om_\ga]_* 
\simeq\HD\II^\om_\ga/[\HD\II^\om_\ga,\,\HD\II^\om_\ga]_*$$ and that 
$$\MM_{q^\half}\cap\II^\com_\ga\simeq\HS\II^\com_\ga/[\HS\II^\com_\ga,\,\HS\II^\com_\ga]_* 
\simeq\HD\II^\com_\ga/[\HD\II^\com_\ga,\,\HD\II^\com_\ga]_*$$  and that
the kernel of the formal $q$-Fourier transform is exactly the commutator
ideal. Hence it makes sense to look for an inverse of the formal $q$-Fourier
transform $\tilde\FF_\ga$ on the image of  $\MM_{q^\half}\cap\II^\om_\ga$ and
$\MM_{q^\half}\cap\II^\com_\ga$.  Moreover, $\tilde F_\ga(f)$ is also a
well-defined function on a neighbourhood of $0$ for
$f\in\MM_s\cap\II^{\rho,s}_\ga$ and $s<q^{-\half}$  and its radius of
convergence  will be $s^{-1}q^{-\half}>1$. $\tilde F_\ga$ is injective on
$\MM_{q^{\half}}\cap\II^{\rho,q^{\half}}_\ga$ by Lemma \ref{tom-odd}. Since
$\tilde\FF_\ga$ is $\mu_\ga$ up to a  multiplicative shift of the variable, 
its inverse boils down to retrieving back a function knowing its $q$-moments. 
This can  clearly be achieved by means of the functions $G_{r,\ga}$'s
as it was
shown  in Section \ref{Gicappa}.\\ Define the operator  $\GG_\ga$ on the
space of functions in $\HD$ whose radius of convergence is greater than
$q^{-1}$ as follows. 
\begin{equation}\label{giggi}\GG_\ga(f)=\GG_\ga\biggl(\sum_{k=0}^\iy
c_kx^k\biggr)=\sum_{k=0}^\iy i^kc_k (q;q)_k
G_{k,\ga}=\Bigl(\sum_{k=0}^\iy(-i)^kc_k(1-q)^k\pa^k
\Bigr)\,u_\ga\end{equation} where $\sum_{k=0}^\iy c_kx^k$ is the power series
expansion of $f$ on a neighbourhood of $0$. By definition  of the
$G_{k,\ga}$'s it is clear that $c_q(\ga)\GG_\ga$ is independent of
$\ga$.\\ Let us
define the following spaces of functions: for $\al>0$
$$\EE_\al^\om:=\{f\in\EE\,|\,f=\sum_kc_kx^k{\hbox { and }} \exists b>0\,|,
|c_k|=O(q^{\half\al k^2}b^k) {\hbox { for }} k\to\iy\}$$ i.e. $\EE_\al^\om$ is
the space of functions of $q$-exponential growth of order $\al$ and
finite type. Let $$\EE^\om:=\bigcup_{\al>0}\EE_\al^\om$$
It is almost tautological that $\tilde\FF_\ga\bigl(\II^{\EE}_{\ga}\bigr)= \EE$,
$\tilde\FF_\ga\bigl(\II^{\rho,s}_{\ga}\bigr)= \HD_{s^{-1}q^{-\half}}$ where
the lower index by $\HD$ will  denote from now on the lower bound of the radius
of convergence (note that in \cite{giotom}  $\HD_a$ meant that
 the radius of convergence had to be greater or equal to $a$ while 
here it denotes that the radius of convergence has to be {\em strictly
greater} than $a$). We also have,
$\tilde\FF_\ga\bigl(\II^{\om}_{\ga,\al}\bigr)= \EE_\al^\om$ and
$\tilde\FF_\ga\bigl(\II^{\om}_{\ga}\bigr)= \EE^\om$.  We get the following
result:  \begin{proposition}\label{funziona} $\tilde\FF_\ga$ defines an
isomorphism of vector spaces between $\MM_{q^\half}\cap\II^{\om}_{\ga,\al}$
and $\EE_\al^\om$ and isomorphisms of algebras between
$\MM_{q^\half}\cap\II^{\om}_{\ga}$ and $\EE^\om$, between
$\MM_{q^\half}\cap\II^{\EE}_{\ga}$ and $\EE$ and between
$\MM_s\cap\II^{\rho,s}_\ga$ and $\HD_{s^{-1}q^{-\half}}$ for $s\le q^\half$.
\end{proposition} \Proof  By the discussion in the previous Sections one finds
that  $\GG_\ga\bigl(\EE)=\II_\ga^\com$, $\GG_\ga\bigl(\EE^\om)=\II_\ga^\om$
and $\GG_\ga\bigl(\HD_{s^{-1}q^{-\half}})=\II_\ga^{\rho,s}$ if $s<1$. By
construction $\tilde\FF_\ga\circ\GG_\ga=\id$ on $\HD_{q^{-1}}$ and
$\GG_\ga\circ\tilde\FF_\ga=\id$  on $\II^{\rho, s}_\ga$ for $s\le q^\half$. 
The rest is clear.\hfill$\square$\smallbreak

\noindent This inversion formula extends the inversion results in
\cite{Koo}. Indeed Koornwinder showed therein that $\FF_\ga$  establishes a particular isomorphism between the
vector space $P_e$ of polynomials times $e_{q^2}(-X^2)$ and the vector space
$P_E$ of polynomials times $E_{q^2}(-q^2X^2)$ resembling the classical case. 
Since $P_e\subset\HS\II^{\som}_\ga$, by Proposition \ref{results2}
$\tilde\FF_\ga=\FF_\ga$ on $P_e$ and since
$P_e\subset\MM_{q^\half}\cap\II^\om_\ga$, the two inverses must coincide on
$P_E$. Koornwinder's inverse transform is essentially given  by 
\begin{equation}(\FF'_\ga f)(y):={1\over{c_q(\ga)\,b_q}}\int_{-1}^1
e_q(ixy)f(x)d_qx\label{toms}\end{equation}  where $c_q(\ga)$ is as in formula
(\ref{eq-moment}). Hence $(\FF'_\ga f)(y)\in\HS$ for every function $f$
bounded in $(-1,\,1)$. I will show that $\FF'_\ga$ and $\GG_\ga$ coincide
on  the whole $\HD_1$. \\ Let $f(x)=\sum_{k=0}^\iy c_k x^k$ for $|x|\le
q^{-1}$. For $|\Im(y)|<1$, by dominated convergence we have
\begin{equation}{1\over{c_q(\ga)\,b_q}}\int_\ga
e_q(ixy)f(x)d_qx={1\over{c_q(\ga)\,b_q}}\sum_{k=0}^\iy
c_k\int_{-1}^1e_q(ixy)x^kd_qx\end{equation} If we denote $q$-differentiation
with respect to $y$ by $\pa_y$ we have for $y\not=0$

\begin{equation}\label{eiqqu}\pa_y\biggl(\int_{-1}^1e_q(ixy)x^kd_qx\biggr)=
{i\over{(1-q)}}\int_{-1}^1e_q(ixy)x^{k+1}d_qx\end{equation} that can be
extended  by continuity at $y=0$. Hence \begin{equation}(\FF'_\ga
f)(y)={1\over{c_q(\ga)\,b_q}}\sum_{k=0}^\iy(-i)^k
c_k(1-q)^k\pa_y^k\int_{-1}^1e_{q}(ixy)d_qx\end{equation} Besides
\begin{eqnarray*}&&\int_{-1}^1e_{q}(ixy)d_qx={{(1-q)}\over{(iy;q)_\iy}}{}_2\phi_1(q,iy;0;q,q)+
{{(1-q)}\over{(-iy;q)_\iy}}{}_2\phi_1(q,-iy;0;q,q)\\
&&={}_2\phi_1(q,0;q^2;q,iy)+
{}_2\phi_1(q,0;q^2;q,-iy)=2\,{}_2\phi_1(0,0;q^3;q^2;-y^2)\\
&&=2\,g_1=2\,c_q(\ga){{(q^2;q^2)_\iy}\over{(q^3;q^2)_\iy}}u_\ga=b_qc_q(\ga)u_\ga
\end{eqnarray*}
where the second  equality follows   by (0.6.13)
 in \cite{KoSw}.
Hence $\GG_\ga$ and $\FF'_\ga$ coincide on $\HD_{q^{-1}}$.  
 As a byproduct we
have found that
$$G_{k,\ga}={{(-i)^k}\over{(q;q)_k}}\FF'_q(x^k)={{(-i)^k}\over{(q;q)_k}}{1\over{b_q\,c_q(\ga)}}\int_{-1}^1x^k e_q(ixy)d_qx.$$
The $G_{k,\ga}$'s are the basis corresponding  to the basis of $\EE$
given by monomials. There holds some sort of orthogonality between
the two bases since  $\int_\ga G_{k,\ga}x^l=q^{-\half(k^2+k)}\delta_{k,l}$.
This explains many of the properties of the $G_{k,\ga}$'s with respect to
product and integration.  Observe also that as a consequence of formula
(\ref{eiqqu}) one can prove by induction that there are polynomials $r_k$ and 
$l_k$ of degree at most $k-1$ for which $y^{k+1}c_q(\ga)G_{k,\ga}=r_k
(y)\cos_q(y)+l_k (y)\sin_q(x)$ where $\cos_q(y):=\half(e_q(iy)+e_q(-iy))$ and
$\sin_q(y):={1\over{2i}}(e_q(iy)-e_q(-iy))$. In some sense then expansion in
terms of  $G_{k,\ga}$ is midway between a $q$-Fourier transform and a
$q$-Fourier series. Interesting results about a $q$-analogue of Fourier series
were obtained in \cite{buslov}, where continuous integrals are
involved. It is interesting that the basic exponentials studied in
\cite{buslov} and references therein depend on two variables, and they
are related to our $g_m$'s for a particular value of the first
variable. The connection between the two families can be the subject
of future research.\\
Observe  that  $\FF'_\ga$ is a
priori an inverse of $\tilde\FF_\ga$ and not of $\FF_\ga$ since
$\MM_{q^\half}\cap\II_\ga^\om\not\subset\II^{\som}_\ga$. On the other hand,
$\tilde\FF_\ga$ and $\FF_\ga$ coincide on the ideal  $I_e$ of
$\MM_{q^\half}\cap\II^\om_\ga$  generated by $e_{q^2}(-X^2)$ because this ideal
is contained in $\HS\II^\som_\ga$ by statement {\sl 2} of Proposition
\ref{results}.  By arguments similar to those in the proof of Proposition 5.3
in \cite{giotom} one shows that the ideal $I'_e$ generated by $e_{q^2}(-X^2)$
in $\MM_{q^\half}\cap\II_\ga^{\EE}$ is contained in
$\MM_{q^\half}\cap\II_\ga^{{s\EE}}$. $I_e$ can be described
explicitely as the space of functions $F=f\,e_{q^2}(-X^2)$ where
$f(x)=\sum _{k=0}^\iy
c_k{\tilde h}_k(x;q)$ for which there are a $c$ and an $\al>0$ such that 
$|c_k|\le C q^{\half(\alpha+1)k^2}c^k$. This is shown using formula $(6.2)$ in
\cite{giotom}. One shows that such an $F$ has to be $F=g*_\ga e_{q^2}(-X^2)$
for some $g\in \HD\II^{\om}_{\ga,\alpha}$, and using the results in the
previous section one sees that $g$ can be chosen to be in
$\MM_{q^\half}\cap\II^\om_{\ga,\al}$. Similarly $I'_e$ can be described
explicitely as the space of functions of the form $f\,e_{q^2}(-X^2)$ where the
coefficients $c_k$ of the power series expansion of $f$ are
$O(q^{\half(k^)}R^k)$ for every $R$ as $k\to\iy$. In the terminology of
\cite{ramis}, this means that $f$ is $q$-Gevrey-Beurling of order $-1$. \\
We have:
\begin{equation}\FF_\ga(I_e)=\tilde\FF_\ga\bigl((\MM_{q^\half}\cap\II^\om_\ga)*_\ga
e_{q^2}(-X^2))=\EE^\om\,E_{q^2}(-q^2X^2)\end{equation} and 
\begin{equation}\FF_\ga(I'_e)=\tilde\FF_\ga\bigl((\MM_{q^\half}\cap\II^\EE_\ga)*_\ga
e_{q^2}(-X^2))=\EE\,E_{q^2}(-q^2X^2)\end{equation} by Proposition
\ref{results2} and the results in \cite{Koo}. In particular, this proves that
$\FF'_\ga(\EE\,E_{q^2}(-q^2X^2))\subseteq
\MM_{q^\half}\cap\II^{s\com}_{\ga,\half}$ and
$\FF'_\ga(\EE^{\om}\,E_{q^2}(-q^2X^2))\subseteq
\MM_{q^\half}\cap\II^{\som}_{\ga,\half}$. On the space  $\EE\,E_{q^2}(-q^2X^2)$
both $\FF'_\ga$ and $\GG_\ga$ coincide with the case $n=1$ and $\ga=1$ of
$F''(\id,\ga)$ in \cite{Ca2} and \cite{giotesi}, with $q^2$ replaced by $q$.
There, the integral is unbounded, but it coincides with a bounded one since
$E_{q^2}(-q^2x^2)=0$ for $x=\pm q^{-k}$ with $k\in{\bf Z}_{\ge1}$. The 
context of  \cite{Ca2} and \cite{giotesi} was the braided setting (see also
\cite{KeMa}) and the integral is also slightly different though. In our
context,  $F''(\id,\ga)$  for $\ga=q$ is the operator $$(\tilde\FF'_q
f)(y)=\sum_{e=0}^\iy{{i^ey^e}\over{(q;q)_e}}\int_q f X^e$$ mapping
$\II^\om_{q,>1}$ to $\EE$ and $\II^\om_{q,1}$ to $\HD$. \\    %
\begin{remark}\label{fourier-remark}\rm If we view $\CC[[x]]$ as the  braided
line  $\cove$ with  the braided Hopf algebra structure recalled in Remark
\ref{braidedline}, then  $F''(\id,\ga)$ is (up to a small change) the braided
Fourier transform defined in \cite{KeMa} but with a nonbosonic  integral. In
this context the variable $y$  lives in $\vece$, another braided Hopf algebra.
$\vece$ is isomorphic, as an algebra, to $\CC[[y]]$, it acts on $\cove$ by
letting $y$ act on $f\in\cove$ as $\pa$ and it has a nontrivial braided Hopf
algebra pairing with $\cove$, hence it is dual to $\cove$. \\   The braided
Hopf algebra structure on $\vece$  is given by the braiding $\Psi(y^k\otimes
y^l)=q^{kl}y^l\otimes y^k$,  the comultiplication
$\Delta(y)=y\otimes1+1\otimes y$, the braided antipode
$S(y^k)=(-1)^kq^{k\choose2}y^k$ and the counit
$\varepsilon(y^k)=\delta_{k,0}$. Since the braiding between $\cove$ and
$\vece$ is nontrivial, $1\otimes y$ and $x\otimes 1$ do not commute in
$\cove\otimes\vece$. In particular  $\psi(x^k\otimes\pa^l)=q^{-kl}\pa^l\otimes
x^k$ and $(x\otimes \pa)^r=q^{-{r\choose2}}x^r\otimes\pa^r$. Hence, formally 
$$\tilde\FF_\ga f=\Bigl(\int_\ga\otimes\id\Bigr)(f\,E_{q}(ix\otimes y))\quad
{\hbox{and }}\quad \tilde\FF'_\ga f=\Biggl(\int_\ga\otimes
SQ\Biggr)(f\,E_q(ix\otimes y))$$  where $Q$ acts on $\vece$ as $q$-shift of
the indeterminate $y$.\hfill$\spadesuit$\end{remark} I  conclude this section
with a rigorous proof of a result in \cite{KeMa} concerning the behaviour of
$\tilde\FF'_\ga$ with respect to the $q$-convolution.

\begin{definition} Let ${\cal C}[[x,\,y]]$ be the space of power series in $x$
and $y$  that converge in some polydisk
$\{x\,|\,|x|<r_1\}\times\{y\,|\,|y|<r_2\}$. We define the operator
\begin{eqnarray*}\Psi\colon{\cal C}[[x,\,y]]&\to& {\cal C}[[x,\,y]]\\
\sum_{n,m}a_{nm}x^n\,y^m&\mapsto& \sum_{n,m} a_{n,m}q^{nm}
y^m\, x^n
\end{eqnarray*} 
\end{definition} \begin{remark}\rm The operator $\Psi$ is, in the
braided context, the braiding from $\cove\otimes\cove$ or $\vece\otimes\vece$
to itself, see also Remark \ref{braidedline}.\hfill$\spadesuit$\end{remark}

\begin{proposition} Let $f\in\HD\II_{\ga,\al}^\om$ with  $\alpha>1$ and let
$g\in\HD\II_{\ga',\beta}^\om$ with $\beta>1$.  If there holds
$(\alpha-1)(\beta-1)>1$ then
\begin{equation}\label{senza-antipodo}\tilde\FF'_{\ga'}(f*_{\gamma}g)=
m\circ\Psi^{-1}(\tilde\FF_{\ga}(f)\otimes\tilde\FF'_{\ga'}(g))\end{equation} 
and it  converges absolutely everywhere.  Equality
(\ref{senza-antipodo}) holds as equality of power series in $\HD$ if 
$(\alpha-1)(\beta-1)=1$ or if  $f$ is as above and $g\in\II_{\ga',1}^\som$.
\end{proposition} {\bf Proof:} Observe that since
$\HD\otimes\HD$ can be embedded in ${\cal C}[[x,\,y]]$,
$\Psi$ is well-defined  on $\HD\otimes\HD$. \\
One sees that $\tilde\FF_\ga=Q\circ
S\circ \tilde\FF_\ga'$, with $S$ as in Definition \ref{antipode}.  Observe
that  $S\circ m=m\circ(S\otimes S)\circ\Psi$ where $m$ is the ordinary
product, $Q$ commutes with $\Psi$ and $S$ on power series and  $Q\circ
m=m\circ (Q\otimes Q)$. By Proposition \ref{results2} 
\begin{eqnarray*}&&S\tilde\FF'_{\ga'}(f*_{\gamma}g)\\  &&=m(S\otimes
S)(\tilde\FF'_{\ga}(f)\otimes\tilde\FF'_{\ga}(g))=S\circ m\circ
\Psi^{-1}(\tilde\FF'_{\ga}(f)\otimes\tilde\FF'_{\ga'}(g)) \end{eqnarray*}
where the composite $S\circ m\circ \Psi^{-1}$ is a well-defined operator on
${\cal C}[[x,\,y]]$ although $\psi^{-1}$ itself may not.\\  If
$(\alpha-1)(\beta-1)>1$ by  Proposition 4.6 in \cite{giotom}  
  $\tilde\FF'_{\ga'}(f*_{\gamma}g)\in{\cal E}$. In this case  we
can multiplyboth sides  by the inverse of the antipode $S$  without problems, obtaining
equality (\ref{senza-antipodo}). Both sides of the equality will be absolutely
convergent everywhere.  If $(\alpha-1)(\beta-1)=1$ or if  $\alpha>1$ and
$g\in\II_{\ga',1}^\som$, then  $\tilde\FF'_{\ga'}(f*_{\gamma}g)$ converges
absolutely on a neighbourhood of zero by statement {\sl 3} of Proposition
\ref{results} and we get the statement. \hfill$\square$\smallbreak
\begin{remark}\rm One could also check that the series
$m\Psi^{-1}(\tilde\FF'_{\ga}(f)\otimes\tilde\FF'_{\ga'}(g))$ converges in a
neighbourhood of zero by direct computation using  the fact that if  the power
series $F(x)=\sum_{n=0}^\infty a_nx^n$ is such that  $|a_n|\le C a^n
q^{{{(\alpha-1)}\over2}n^2}$ for some $C,\,a\ge0$ and some $\alpha>1$ and if
the power series $G(x)=\sum_{n=0}^\infty b_nx^n$ is such that  $|b_n|\le B b^n
q^{{{(\beta-1)}\over2}n^2}$ for some $B,\,b\ge0$ and some $\beta>1$ and if
$(\alpha-1)(\beta-1)>1$ then $$m\circ\Psi^{-1}(F(x)\otimes
G(x))=\sum_{t=0}^\infty\biggl[\sum_{m+n=t}a_m b_n
q^{-nm}\biggr]x^t\quad{\hbox{and}}$$ 
\begin{eqnarray*}&&\biggl|\sum_{m+n=t}a_m
b_n q^{-nm}\biggr|\\
&&\le CB \sum_{m+n=t}a^mb^n
q^{\half(
(\alpha-1)-{1\over{(\beta-1)}})m^2}
q^{{1\over{2(\beta-1)}}(m^2+(\beta-1)^2n^2-2(\beta-1)nm)}\\
&&\le D (\max
(a,\,b))^t
\sum_{m=0}^tq^{\half(
(\alpha-1)-{1\over{(\beta-1)}})m^2}\\
\end{eqnarray*}\hfill$\spadesuit$\end{remark}
\begin{remark} \rm In  \cite{Zhang}
, where $q>1$,  an analytic inversion
formula along a direction in $\bf C$ for the $q^{-1}$-Borel transform is
given  for functions with a particular growth.  Recall that on $\HD$, $S$ is
essentially the $q$-Borel transform. See also \cite{ramis} for the
isomorphism of vector spaces between $\HD$ and $q$-Gevrey series.

\hfill$\spadesuit$\end{remark}

 \section{Invertibility of
functions}\label{invertsection}

Next question is whether given $f\in\HD\II^{\EE}_\ga$  there exists a
function for which $f*_\ga g=u_\ga$ and in that case whether $g*_\ga f$ is also
well-defined and  equal to $u_\ga$.  If a left  inverse exists, it will
not be unique since any element in $g+[\HD\II^\com_\ga,\HD\II^\com_\ga]_*$ 
will also be a left inverse.\\ By formula (\ref{generating}) it follows that a
necessary condition for invertibility of  a function of left type is that
$\mu_{0,\ga}(f)=\int_\ga f\not=0$. In particular, odd functions can never be
invertible. \\ Observe also that if  an inverse of $\mu_\ga(f)$ exists, then
it might not correspond to a function of  $\II^\com_\ga$, since the only
functions that are invertible in $\cal E$ are those with no zeroes. \\ On
the other hand  if $\mu_\ga(f)(t)\not=0$ for $|t|<\rho$ then the
inverse of $\mu_\ga(f)$ will be defined and analytic for $|t|<\rho$ and  its
power series expansion on this disk will be \begin{equation}
\nu(t)=\sum_{k=0}^\iy{{d_kt^k}\over{[k]_q!}}\quad {\hbox { and
$|d_k|=O({\sigma^k})$ for $k\to\iy$ for every
$\sigma>\rho^{-1}(1-q)^{-1}$.}}\end{equation}  In Section \ref{Gicappa} we
showed how to construct  a function
$g\in\II_\ga^\iy$ such that $\mu_\ga(g)(t)=\nu(t)$ on a neighbourhood of $0$
at least if $\rho$ is big enough. \\  Observe
that if $f\in\II_\ga^\com$ then for any $p\in\ZZ$ we have 
$Q^{-p}f\in\II_\ga^\com$  and $\mu_\ga(Q^{-p}f)(t)=q^{p}\mu_\ga(f)(q^{p}t)$,
because $\int_\ga Q(F)=q^{-1}\int_\ga F$ for every $F\in\II_\ga$.  Hence if
$\mu_\ga(f)(t)\not=0$ for $|t|<\rho$, then $\mu_\ga(Q^{-p}f)(t)\not=0$ for
$|t|<\rho q^{-p}$. This tells that the conditions on $f$ are satisfied at
least for its (big) $q$-shifts. \begin{proposition}\label{inversion} Let
$f\in\II_\ga^\om$ and let $\mu_\ga(f)(t)\not=0$ for $|t|<\rho$. If
$\rho>(1-q)^{-1}$ there exists a function $g\in\MM_{s}$ with $s<q^{-\half}$
such that $\mu_\ga(g)\mu_\ga(f)=1$.  If  moreover, $\rho>q^{-1}(1-q)^{-1}$ then
$g\in\MM_{q^\half}\cap\II^{\rho,q^{-1}}_\ga$. In this case, $f*_\ga
g=u_\ga$.\end{proposition} \Proof  It is a consequence of Lemma's \ref{appa},
\ref{construct2} and \ref{construct3}. The fact that
$f*_\ga g=u_\ga$ follows from Lemma \ref{tom-odd}.\hfill$\square$\smallbreak 
\begin{corollary} Let $f$, $g$, $\rho$ be defined as in Proposition
\ref{inversion} with $\rho>(1-q)^{-1}q^{-1}$ and  let $f\in\MM_{q^\half}$. Then
$f*_\ga g=g*_\ga f=u_\ga$. \end{corollary}
\Proof  We  have to show that $g*_\ga f=u_\ga$. By
Lemma \ref{construct3} $g*_\ga f$ is a well-defined function
in $\MM_{q^\half}$, hence $g*_\ga f\in\II_\ga^\iy$. Formula
(\ref{generating}) holds by Corollary \ref{new-algebra}, so that by
Lemma \ref{tom-odd}  the statement follows. \hfill$\square$ \smallbreak 
\begin{proposition}Let $g\in\II_\ga^\iy$ be such
that $\mu_\ga(g)=\sum_{k=0}^\iy{{d_k}\over{[k]_q!}}t^k$ with
$|d_k|=O(\sigma^k)$ for $k\to\iy$, and let $h$ be a function defined,
together with its $q$-derivatives, on a domain $\Omega$. If there is an
$M\in(0,\,\sigma^{-1}(1-q)^{-1})$ for which $|\pa^kh(x)|=O(M^k)$ for
$x\in\Omega$,  then $g*_\ga h$ is well-defined on $\Omega$. In
particular, if $h\in\HD_R$ for $R>\sigma$ then
$g*_\ga h\in\HD_R$, and if $h\in{\cal E}$, then $g*_\ga
h\in{\cal E}$.\end{proposition} \Proof  The proof
follows like the proof of Lemma 3.4 in \cite{giotom}.\hfill$\square$\smallbreak

\begin{corollary}\label{inverse}Let $f$, $g$ and $\rho$ be defined as in
Proposition \ref{inversion}, with $\rho>(1-q)^{-1}q^{-1}$. Let $h\in\HD_R$ for
  some $R>\rho^{-1}(1-q)^{-1}$. Then $g*_\ga h\in\HD_R$. 
In particular, if $h\in{\cal E}$  then $g*_\ga h\in{\cal E}$. \end{corollary}
\Proof  Clear by the previous results.\hfill$\square$\smallbreak
\section{Convolution and $q$-differential equations}\label{equadiffe}

In this section I will  show a link between invertibility of a function in
$\MM_{q^\half}\cap\II^\com_\ga$ and the solution of a $q$-differential equation
  with constant coefficients on a fixed
$q$-lattice $L(\ga)$. A {\em $q$-differential equation} is  an equation of
the form $L\,Y=F$ where $F$ is a given function , $Y$ is an unknown and
$L\in{\bf C}[\pa]\subset{\bf C}[x^{-1},\,Q]$. Hence such an equation can be
reduced to a $q$-difference equation with polynomial coefficients and with
leading term $=1$  by multiplying both sides by
$(-1)^nc_n^{-1}(1-q)^nq^{n\choose2}x^n$ where $n$ is the order of the
equation and $c_n$ is the coefficient of $\pa^n$ in $L$. The associated
$q$-difference equation are always  regular singular at $0$ (see \cite{ramis}
for a definition) and their characteristic  equation (see \cite{adams3}) has
roots $1,\,q,\,\ldots,\,q^{n-1}$. Many interesting results about the behaviour
of the solutions were known already in the 30's, see \cite{adams3} and
\cite{adams2}, where in particular the solutions for homogeneous equations
are described.  Many methods for solving $q$-difference equations are known,
hence I  do not think that using $q$-convolution shall simplify the  problem
in general. It  offers though a different interpretation and can help
to describe the type of solutions one could find and to find a
particular solution in special cases, as I  shall show.\\  Let  
\begin{equation}\label{prima-eq}L\,Y=\sum_{n=0}^Nc_k\pa^k Y=F\end{equation}
for a given function $F$ and an unknown $Y$. Equation (\ref{prima-eq}) is
equivalent to \begin{equation}\label{seconda-eq}\sum_{n=0}^Nc_k\pa^ku_\ga *_\ga
Y=\Biggl(\sum_{k=0}^N(-1)^kc_k[k]_q! G_{k,\ga}\Biggr) *_\ga Y= D_L*_\ga
Y=F.\end{equation}  where $D_L$ is of left type for every $\al>0$ since
$\mu_{\ga}(D)$ is a polynomial.\\ If $c_0\not=0$ we can invert 
$\mu_\ga(D_L)=\sum_{n=0}^N(-1)^nc_nt^n$ on a neighbourhood of $0$.\\ If
$c_0=c_1=\cdots=c_{l-1}=0$ and $c_{l}\not=0$, then equation (\ref{prima-eq})
becomes \begin{equation}\label{terza-eq}\sum_{p=0}^{N-l}c_{p+l}\pa^p (\pa^l
Y)=F.\end{equation}  If we can solve $\sum_{p=0}^{N-p}c_{p+l}\pa^p (Z)=F$, 
and if the solution behaves well (for instance, if it belongs to $\HD$), then
we can determine solutions of equation (\ref{terza-eq}) applying $l$ times 
indefinite  $q$-integration  $\int_0^x f(t)d_qt$. As in classical integration,
 $q$-integration determines $q$-primitives up to a constant.\\  Now suppose
$c_0\not=0$. $\mu_\ga(D_L)$ is a polynomial, hence it will have zeroes. If
$\mu_\ga(D_L)(t)\not=0$ for $|t|<\rho$ and $\rho>(1-q)^{-1}$, then we can
construct $g$, a left inverse of $D_L\in\MM_{q^{\half}}\cap\II^\om_\ga$ and it
will belong to $\MM_{s}$ for some $s<q^{-\half}$.  The function $g$ should be
seen as a  $q$-analogue of  the classical fundamental solution   associated to
the differential equation $Ly=F$, since for a fundamental solution there should
hold: $Lg=\delta$ where $\delta$ is Dirac's delta.  By the results in the
preceding sections, if $\rho$ is big enough we may compute $g*_\ga F$. Even if
the function $g$ does not belong to $\II^\iy_\ga$,
we still  can formally compute $g*_\ga F$ using the coefficients of its
expansion with respect to the $G_{k,\ga}$'s as if they were really $q$-moments.
Then  \begin{itemize} \item If $F\in\MM_s$ with $s<q^{-\half}$, then  the
formal product $g*_\ga F\in\MM_{r}$ for  some $r<q^{-\half}$. 
\item If $F\in
\MM_s$ for $s<q^{-\half}$ and $\rho>(1-q)^{-1}s^{-1}q^{-\half}$ then $g*_\ga
F\in\MM_{s}$.
 \item If $F\in\HD_{\rho'}$ for $\rho'>\rho^{-1}(1-q)^{-1}$, then
$g\in\II_\ga^\iy$ and  $g*_\ga F\in\HD_{\rho'}$. 
\item If $F\in{\cal E}$ then
$g*_\ga F\in{\cal E}$. \end{itemize} In all those cases,  both $D*_\ga(g*_\ga
F)$ and $(D*_\ga g)*_\ga F$ are well-defined, hence they are both equal to $F$
by dominated convergence. This implies that  $g*_\ga F$ is a solution of
equation (\ref{prima-eq}). If we replaced $g$ by another function $g'$ for
which $\mu_\ga(g')\mu_\ga(f)=1$ on a neighbourhood of zero, then $g*_\ga
h=g'*_\ga h$ for every $h$ for which the product is defined, hence we would
get the same solution. Moreover, unicity of a possible solution in $\MM_s$
for $s<1$ follows by determinacy of the $q$-moment problem on $\MM_s$. If $F$
is a polynomial instead, the solution obtained by this construction will be
again a polynomial. \\ Suppose now that $c_0\not=0$ but  $\rho\le
(1-q)^{-1}q^{-1}$. We can find a $p\in\Zplus$ for which 
$\rho_p:=q^{-p}\rho>(1-q)^{-1}q^{-1}$. Observe that $q^{-p}\rho$ is related 
to the homogeneous equation
\begin{equation}\label{quarta-eq}\Bigl(\sum_{n=0}^Nc_nq^{pn}\pa^n\Bigr)
Y=0\end{equation}  as $\rho$ is related to equation (\ref{prima-eq}). In this
case, the role of $D$ is played by  the function
$D_p:=\sum_{k=0}^N(-1)^kc_k[k]_q!
q^{pk}G_{k,\ga}\in\MM_{q^\half}\cap\II_\ga^\om$. The inverse of
$\mu_\ga(D_p)(t)=\mu_\ga(D)(q^pt)$  is analytic for $|t|<q^{-p}\rho=\rho_p$.\\
The corresponding  function  $g_p$ is well-defined and it is such that
$D_p*_\ga g_p=u_\ga$ by Corollary \ref{inverse}. If  $F\in\HD_{\rho'}$ with
$\rho'>\rho^{-1}(1-q)^{-1}$, then the function $F_p:=Q^{-p}F\in\HD_{\rho'_p}$
where $\rho'_p=\rho q^{p}>q^p\rho^{-1}(1-q)^{-1}=\rho_p^{-1}(1-q)^{-1}$. Hence
  $Y_p:=g_p*_\ga F_p\in\HD_{\rho'_p}$ and it is a solution of:
\begin{equation}\label{quinta-eq}\sum_{n=0}^Nc_nq^{pn}\pa^n
Y=F_p.\end{equation} Then $Q^pY_p\in\HD_{\rho'}$ and 
\begin{equation}\sum_{n=0}^Nc_n\pa^nQ^{p}Y_p=\sum_{n=0}^Nc_n
q^{np}Q^p\pa^nY_p=Q^pF_p=F\end{equation} so that $Q^pY_p$ is a solution of
equation (\ref{prima-eq}). This shows that if $\rho'$ is big enough, the
inverse of $D$ may not be well-defined but we could still apply this method in
order to find a particular solution.    

\begin{example}\rm Consider the
equation
\begin{equation}\label{esempio}Y-q^{2r}\pa^2Y=e_{q^2}(-x^2).\end{equation}
This equation is equivalent to $D*_\ga Y=e_{q^2}(-x^2)$ with
$D=u_\ga-q^{2r}[2]_q!G_{2,\ga}$. \\ $\mu_\ga(D)(t)=1-q^{2r}t^2$, so that its
inverse for $|t|<q^{-2r}$ is equal to $\nu(t)=\sum_{p=0}^\iy
(q^{2r}t^2)^p=\sum_{p=0}^\iy q^{2rp}t^{2p}$. For $r$ big enough,
$\rho=q^{-r}>(1-q)^{-1}$, hence $$g:=\sum_{p=0}^\iy
[2p]_q!q^{2rp}G_{2p,\ga}=\sum_{p=0}^\iy q^{2rp}\pa^{2p}u_\ga$$ is
well-defined, as well as  $$g*_\ga e_{q^2}(-x^2)=\sum_{p=0}^\iy
q^{2rp}\pa^{2p}e_{q^2}(-x^2)=\sum_{p=0}^\iy
{{q^{2rp}q^{2p^2-p}}\over{(1-q)^{2p}}}{\tilde h}_{2p}(x;q)e_{q^2}(-x^2)$$
which is a solution of the equation. \hfill$\spadesuit$\end{example}  If
$F\in\MM_s$ with $s<1$ then the method of $q$-shift does not apply since
$\rho'=1$ can never be greater than  $\rho^{-1}(1-q)^{-1}>1$.  In this case,
one can  approximate the solution by multiplying $F$ by finite linear
combinations $g_n$ of $G_{k,\ga}$'s converging to $g$ in the topology
described in the previous Section.  Then since the $q$-convolution product is
continuous with respect to this topology,  $g_n*_\ga F$ will converge to $Y$
with respect to this topology. \\

\begin{remark}\rm The approach we used to solve a  $q$-differential
  equation 
 can be translated into applying $\tilde\FF_\ga$ to
both sides of  equation (\ref{prima-eq}) and considering whether $\GG_\ga$ can
be applied to
$\bigl(\tilde\FF_\ga(D)\bigr)^{-1}\tilde\FF_\ga(F)$.\hfill$\spadesuit$\end{remark}

\section{Acknowledgements} The
author wants to thank T. Koornwinder for many stimulating discussions, for his
encouragement  and for interesting suggestions. In particular, Lemma \ref{tom}
and  Section \ref{distributions}  are deeply based on his ideas.\\ The
author also wishes to thank the University of Trieste for the
financial support and the Department of Mathematics of the University
of Cergy-Pontoise for the hospitality.

\bibliographystyle{alpha}
\bibliography{libraz}

\begin{thebibliography}{Car99b}

\bibitem[Ada25]{adams1}
C.~R. Adams.
\newblock {Note on the existence of analytic solutions of non homogeneous
  linear $q$-difference equations, ordinary and partial}.
\newblock {\em Annals of Mathematics}, 25:73--83, 1925.

\bibitem[Ada29]{adams2}
C.~R. Adams.
\newblock {On the linear ordinary $q$-difference equation}.
\newblock {\em Annals of Mathematics}, 30:195--205, 1929.

\bibitem[Ada31]{adams3}
C.~R. Adams.
\newblock {Linear $q$-Difference Equations}.
\newblock {\em Bull. AMS}, 31:361--382, 1931.

\bibitem[B{\'e}z93]{bezi}
J.~P. B{\'e}zivin.
\newblock {Sur les \'equations fonctionelles aux $q$-differences}.
\newblock {\em Aequationes Mathematicae}, 43:159--176, 1993.

\bibitem[BS98]{buslov}
J.~Bustoz and S.~K. Suslov.
\newblock {Basic analog of Fourier series on a $q$-quadratic grid}.
\newblock {\em Methods and Applications of Analysis}, 5:1--38, 1998.

\bibitem[Car99a]{giotesi}
G.~Carnovale.
\newblock {\em {Algebraic and Analytic Aspects of the Quantum Yang-Baxter
  equation}}.
\newblock 1999.
\newblock PhD Dissertation, University of Utrecht.

\bibitem[Car99b]{Ca2}
G.~Carnovale.
\newblock {On the braided Fourier transform in the $n$-dimensional quantum
  space}.
\newblock {\em J. Math. Phys.}, 40(11):5972--5997, 1999.

\bibitem[CK99]{giotom}
G.~Carnovale and T.~Koornwinder.
\newblock {A $q$-analogue of convolution on the line}.
\newblock Technical report, 1999.
\newblock math. CA/9909025.

\bibitem[GR90]{gasper}
G.~Gasper and M.~Rahman.
\newblock {\em Basic hypergeometric series}.
\newblock Cambridge University Press, Cambridge, 1990.

\bibitem[KM94]{KeMa}
A.~Kempf and S.~Majid.
\newblock {Algebraic $q$-integration and Fourier theory on quantum and braided
  spaces}.
\newblock {\em J. Math. Phys.}, 35(12):6802--6837, 1994.

\bibitem[Koo90]{Koor}
T.~Koornwinder.
\newblock {Orthogonal polynomials in connection with quantum groups}.
\newblock In P.~Nevai, editor, {\em {Proceedings NATO ASI on Orthogonal
  polynomials, Columbus, Ohio, 1989}}, pages 257--292. Kluwer Academic Press,
  1990.

\bibitem[Koo97]{Koo}
T.~Koornwinder.
\newblock {Special Functions and $q$-commuting variables}.
\newblock In M.E.H Ismail, D.R. Masson, and M.~Rahman, editors, {\em {Special
  Functions, $q$-series and related topics - Fields Institute Communications
  14}}, pages 131--166. AMS, 1997.

\bibitem[Koo99]{informal}
T.~Koornwinder.
\newblock {Some simple applications and variants of the $q$-binomial formula}.
\newblock 1999.
\newblock informal paper.

\bibitem[KS92]{KooSw}
T.~Koornwinder and R.~F. Swarttouw.
\newblock {On $q$-analogues of the Fourier and Hankel transforms}.
\newblock {\em Trans. Am. Math. Soc.}, 333:445--461, 1992.

\bibitem[KS98]{KoSw}
R.~Koekoek and R.~Swarttouw.
\newblock {The Askey-scheme of hypergeometric orthogonal polynomials and its
  $q$-analogue}.
\newblock Technical report, Delft University of Technology, Faculty TWI, 1998.
\newblock Report 98-17.

\bibitem[Maj93]{Majiprimo}
S.~Majid.
\newblock {Braided momentum in the $q$-Poincar\'e group}.
\newblock {\em J. Math. Phys.}, 34:2045--2058, 1993.

\bibitem[Maj95]{Maj1}
S.~Majid.
\newblock {\em {Foundations of quantum groups}}.
\newblock Cambridge University Press, 1995.

\bibitem[OR97]{OR}
M.~Olshanetsky and V.~Rogov.
\newblock {The $q$-Fourier transform of $q$-distributions}.
\newblock Technical report, IHES, December 1997.
\newblock QA/9712055.

\bibitem[Ram92]{ramis}
J.~P. Ramis.
\newblock {About the growth of entire functions solutions of linear algebraic
  $q$-difference equations}.
\newblock {\em Ann. de la Fac. de Toulouse, S\'erie 6}, I:53--94, 1992.

\bibitem[Ryd21]{ryde}
F.~Ryde.
\newblock {\em {A contribution to the theory of linear homogeneous geometric
  difference equations ($q$-difference equations )}}.
\newblock 1921.
\newblock PhD Dissertation, Lund.

\bibitem[Zha99]{Zhang}
C.~Zhang.
\newblock {D\'eveloppements asymptotiques $q$Gevrey et s\'eries
  G$q$-sommables}.
\newblock {\em Ann. Inst. Fourier}, 49:227--261, 1999.

\end{thebibliography}

\end{document}